\documentclass[a4,11pt]{article}

\usepackage{graphicx}
\usepackage{mathtools}
\usepackage{dsfont}
\usepackage{amsmath, amsthm, amssymb}
\usepackage{authblk}

\title{On the Functoriality of Belief Propagation Algorithms on Finite Partially Ordered Sets}

\author[1,2]{Grégoire Sergeant-Perthuis\thanks{\texttt{gregoire.sergeant-perthuis@sorbonne-universite.fr}}}
\author[3]{Toby St Clere Smithe}
\author[1]{Léo Boitel}

\affil[1]{CQSB, Sorbonne Université, Paris, France}
\affil[2]{Ouragan, Inria Paris, France}
\affil[3]{VERSES Research, Topos Institute}

    \theoremstyle{plain}
 \newtheorem{defn}{Definition}
  
  \newtheorem*{defn*}{Definition}
 
  \theoremstyle{plain}
  \newtheorem{thm}{Theorem}
  
  \newtheorem{prop}{Proposition}
  \newtheorem*{prop*}{Proposition}
  \newtheorem{lem}{Lemma}
   \newtheorem*{cor*}{Corollary}
  \newtheorem*{theo*}{Theorem}
  \newtheorem*{thm*}{Theorem}
  
  \theoremstyle{remark}
  \newtheorem{rem}{Remark}
 
 \usepackage[toc,page]{appendix}
\usepackage{calrsfs}
\usepackage{tikz-cd}

\newcommand{\p}{\mathbb{P}}

\newcommand{\R}{\mathbb{R}}
\newcommand{\A}{\mathcal{A}}

\usepackage{pstricks,pst-node,pst-tree}

\usepackage{bm}

\usepackage{relsize}

\newcommand\independent{\protect\mathpalette{\protect\independenT}{\perp}}
\def\independenT#1#2{\mathrel{\rlap{$#1#2$}\mkern2mu{#1#2}}}

 \usepackage{xcolor}

\newcommand{\im}{\operatorname{im}}

\theoremstyle{plain}

\numberwithin{equation}{section}

\usepackage{multiaudience}
\SetNewAudience{long}
\SetNewAudience{GSP}

\usepackage{hyperref}

\begin{document}

\maketitle

\begin{abstract}

Undirected graphical models are a widely used class of probabilistic models in machine learning that capture prior knowledge or putative pairwise interactions between variables. Those interactions are encoded in a graph for pairwise interactions; however, generalizations such as factor graphs account for higher-degree interactions using hypergraphs. Inference on such models, which is performed by conditioning on some observed variables, is typically done approximately by optimizing a free energy, which is an instance of variational inference. The Belief Propagation algorithm is a dynamic programming algorithm that finds critical points of that free energy. Recent efforts have been made to unify and extend inference on graphical models and factor graphs to more expressive probabilistic models. A synthesis of these works shows that inference on graphical models, factor graphs, and their generalizations relies on the introduction of presheaves and associated invariants (homology and cohomology groups). We propose to study the impact of the transformation of the presheaves onto the associated message passing algorithms. We show that natural transformations between presheaves associated with graphical models and their generalizations, which can be understood as coherent binning of the set of values of the variables, induce morphisms between associated message-passing algorithms. It is, to our knowledge, the first result on functoriality of the Loopy Belief Propagation.

\end{abstract}

\section{Introduction and related work}
\label{sec:typesetting-summary}

\subsection{Context}
\textit{Undirected graphical models}, also called \textit{Random Markov fields}, are probabilistic models that relate conditional independence relations between random variables and connectivity of associated nodes \cite{Pearl1988} in a given graph $G=(V,E)$. They are still widely used in machine learning today as they achieve the best performance for data with \emph{a priori} knowledge on dependencies between variables or causal relationships, e.g., Hidden Markov models and Bayesian networks \cite{gimpel2008statistical}. Thanks to the celebrated Hammersley-Clifford theorem, graphical models can be reformulated as particular ways of factoring a distribution over all the variables $X_v; v \in V$. This leads to considering the more general notion of a factor graph between variables \cite{wainwright2008graphical,10.5555/1592967} which encode higher order interactions between variable using hypergraphs. A hypergraph is a bipartite graph where one kind of node corresponds to variables and the other corresponds to a subset $a \subseteq I$ of variables, with directed edges given by the membership relationship $a \to i \iff i \in a$. A factor graph is the data of a hypergraph, $\mathcal{H}$, and factors on each subset $f_a: E_a \to \mathbb{R}; a\in \mathcal{H}$ made of functions that only depend on the variables $X_i, i \in a$.

Given observations for a subcollection of variables \( J \subseteq I \), one wants to compute the posterior $\mathbb{P}_{X_J \mid X_{\overline{J}}}(x_J \mid x_{\overline{J}})$. Computing this posterior amounts to computing the likelihood of non-observed (hidden) variables; we will therefore refer to it as \emph{inference}. Inference on undirected graphical models and factor graphs is exponentially costly with respect to the number of variables, \( |I \setminus J| \). It requires computing the marginal likelihood \( \mathbb{P}_{X_{\overline{J}}} \), which involves summing over all possible value configurations of \( X_J \), consisting of \( N^{|J|} \) terms, assuming that all \( X_i \) take values in a set of cardinality \( N \).

One resorts instead to approximate inference, which is inference up to some error bound. It is achieved through variational inference by observing that inference can be expressed as minimizing a free energy \cite{10.5555/3008751.3008848}, called the Bethe free energy, under convex constraints that account for the factor graph. These constraints encode the compatibility of the marginals of the probability distributions over the variables \(X_e\), where \(e \in \mathcal{H}\) is a hyperedge of the hypergraph on their common vertices. The Bethe free energy is built by adding entropies \( S(p_a) \) minus energy terms \( \mathbb{E}_{p_a}[H_a] \) on the marginal distribution \( p_a \) over \( X_a \), with a discount factor \( c(a) \) that accounts for over-counting some subsets \( a \) when considering elements \( b \) that contain \( a \). More explicitly, \( c(a) \) derives from the inclusion-exclusion formula on sets. When the number of variables tends to infinity, both the entropy of the joint free energy of the joint distribution, $S(p_I) - \mathbb{E}_{p_I}[H]$, and the Bethe free energy coincide.

Mathematically, the set of constraints corresponding to a given undirected graphical model or factor graph arises from finding sections of a certain presheaf \cite{PeltrePhD}, called a `graphical presheaf' in this paper. A presheaf has as its source a collection of points and arrows between them—here, a collection of subsets of variables ordered by inclusion—and maps these subsets to sets and the arrows to functions.  Recent work \cite{sergeantperthuis2022regionalized, sergeantperthuis:hal-04527780} makes use of a broader collection of presheaves for probabilistic modeling of partial and heterogeneous signals, where the compatibility relations are given a priori by stochastic maps. Their setting generalizes, in a unified manner, inference over graphical models and factor graphs by introducing a free energy and a message-passing algorithm (\text{MP}) for inference on their presheaves. Observing that undirected graphical models and factor graphs have a topological interpretation, thanks to these presheaves, allows us to apply a wide range of topological transformations to those presheaves, some of which may change the topology of the functor, while others do not (see sheaf cohomology \cite{gallier2022homology}). Building on \cite{yedidia2005, PeltrePhD, sergeantperthuis2022regionalized, sergeantperthuis:hal-04527780}, we ask how those topological transformations impact the associated critical points of the Bethe free energy and, therefore, the inference of associated probabilistic models. 
Previously the study of the impact of the topology of the factor graph on the critical points of the Bethe free energy can be found in \cite{welling2012choice} and \cite{Peltre}. In the first, the authors study how pruning the poset can keep the free energy unchanged. In the second, a more general notion of acyclicity is given for the poset, which allows extending the proof that the free energy has only critical points. In both cases, although the operations studied and the conditions imposed are of a topological nature, they seem to be isolated results focused on very particular aspects of the poset. Methods such as homotopy continuation can be applied to Belief Propagation in order to modify the fixed points of the Belief Propagation algorithm into those of a simpler algorithm in a continuous fashion \cite{grim2024message}. The homotopy continuation is performed in the space of all algorithms: one follows a straight line between two Belief Propagation algorithms, but one cannot expect that the algorithm along this line is itself a Belief Propagation algorithm for some graph or factor graph.

In our approach, we show how natural transformations arise between presheaves associated with the probabilistic models discussed above. Natural transformations are generic transformations of presheaves widely used in their study and can be applied to any presheaf without requiring additional hypotheses such as acyclicity. In our setting, they specify compatible ways of binning values of groups of variables, similar to a "local coarse graining". We show that a natural transformation \(\phi: F \to G\) translates into a map between the associated message-passing algorithms \(MP_F \to MP_G\), mimicking intertwining of linear operators but, in our setting, for the nonlinear operators defined by the message-passing algorithms (Theorem \ref{main:thm1}). Such transposition is called functoriality. 

\subsection{Structure of the paper}
In Section \ref{preliminarie}, we recall the necessary background on undirected graphical models, factor graphs, their generalizations, and the respective inference and Belief Propagation algorithms. In Section \ref{classical} we define graphical models and factor graphs. In Section \ref{presheaves}, we introduce the Bethe free energy, which is central to inference on undirected graphical models and factor graphs. We introduce graphical presheaves and explain how general presheaves from finite partially ordered sets taking values in finite sets serve as a richer class of probabilistic models. We describe the associated Bethe free energy and related presheaves used for defining message-passing algorithms that extend the Belief Propagation algorithm. In Section \ref{belief-propagation}, we recall the expression of the Belief Propagation update rule for undirected graphical models and its generalization. We explain the correspondence between fixed points of the algorithms and critical points of the Bethe free energies. We recall how the Belief Propagation operator decomposes into several elementary operators, which will prove to be central in showing the functoriality of the Belief Propagation algorithm.

In Section \ref{Main-results}, we present the main contributions of the paper. We show that the message-passing operator transforms in an adequate manner under a natural transformation of the underlying presheaf, i.e., that the message-passing operators are functorial, sending functors $F$ to their non-linear inference update rules $MP_F$ and natural transformations $\phi:F\to G$ to particular linear maps on the operators $MP_F\to MP_G$, see Theorems \ref{main:thm1},\ref{fonctoriality-increment},\ref{label:thm-3}.

In the appendices, we recall definitions related to presheaves from a poset to sets and vector spaces (Appendix~\ref{sheaves}), and the associated Möbius inversions (Appendix~\ref{mobius-inversion-functor}). We also discuss the correspondence between the fixed points of the Generalized Belief Propagation and the critical points of the Bethe Free Energy (Appendix~\ref{proof-old}). We then explain the link between the Message Passing algorithm considered in this article and the General Belief Propagation algorithm (Appendix~\ref{comparison-algorithm}) and their fixed points. Finally, we provide a proof of the main theorems of this article in Appendix~\ref{main-proofs}.

\section{Preliminaries}\label{preliminarie}
\subsection{Graphical Models, Hypergraphs, Factor Graphs}\label{classical}

Let $G = (V, E)$ be an undirected graph, where $V$ is a finite set of vertices and $E$ is a collection of undirected edges, i.e., subsets of $V$ of cardinality $2$. For $v \in V$, let us denote by $\partial v$ the neighbors of $v$, i.e., the elements $v_1 \neq v$ such that $\{v_1, v\} \in E$. A subset $a$ of $V$ such that every two distinct vertices in the subset $a$ are neighbors is called a clique; we denote the set of cliques of $G$ by $\mathcal{C}(G)$. Associate to each vertex $v \in V$ a unique random variable $X_v$ taking values in a finite set denoted by $E_v$. For a subset of vertices $a \subseteq V$, we denote by $X_a$ the associated joint variable $(X_v;v \in a)$, by $x_a$ a configuration of $X_a$, and by $E_a$ the associated configuration space for $X_a$, i.e., $E_a = \prod_{v \in a} E_v$. Following these conventions, as an example, $x_a \in E_a$ and $E_V = \prod_{v \in V} E_v$. Let us denote the space of joint probability distributions of all the random variables $X_v$, with $v \in V$, as $\mathbb{P}(E_V)$; $\p_{>0}(E_V)$ will correspond to the space of strictly positive probability distributions, i.e. for any $x\in E_V$, $P(x)>0$. An \emph{undirected graphical model}, also known as a \emph{Random Markov field} \cite{Lauritzen, lauritzen2019lecture,bishop2006pattern}, with respect to $G$, is a joint probability distribution $P \in \mathbb{P}(E_V)$ such that for any vertex $v \in V$, the random variable $X_v$ is conditionally independent of $X_{V \setminus \{v \cup \partial v\}}$ given $X_{\partial v}$, denoted as $X_v \independent X_{V \setminus (v \cup \partial v)} \mid X_{\partial v}$. The Hammersley-Clifford theorem \cite{Lauritzen} asserts that, for a strictly positive joint probability distribution $P\in \p_{>0}(E_V)$, the joint distribution of a undirected graphical model factors on the cliques of the graphs, i.e. that there is a collection of real valued functions $(f_a\in \R^{E_a}_{>0},a\in \mathcal{C}(G))$ such that, $\forall x\in E_V, P(x)=\prod_{a\in \mathcal{C}(G)} f_a(x_a)$.

\textit{Factor graphs} are an extension of undirected graphical models, starting with a \emph{hypergraph} rather than a graph. A hypergraph, denoted as $\mathcal{H} = (V, \mathcal{H}E)$, is defined by a finite set of vertices $V$ and a collection of finite subsets of $V$ called hyperedges \cite{10.5555/2500991}. Undirected Graphs are particular cases of hypergraphs, where each $a \in \mathcal{H}E$ has cardinality 2. Hypergraphs are themselves particular instances of \emph{partially ordered sets}, which are sets equipped with a binary relation $\leq$ that is transitive (i.e., if $c \leq b$ and $b \leq a$, then $c \leq a$) relfexive ($x\leq x$) and antisymmetric (i.e., if $b \leq a$ and $a \leq b$, then $a = b$); for hypergraphs, the order relation is given by the belonging relation: $v \leq a$ whenever $v \in a$, with $v \in V$ and $a \in \mathcal{H}E$. We denote the poset associated with a hypergraph as $\mathcal{A}(\mathcal{H})$, and the poset associated with a graph, seen as a hypergraph, as $\mathcal{A}(G)$. A factor graph \cite{wainwright2008graphical,10.5555/1592967} is a hypergraph with a unique set of variables $X_v$, for $v \in V$, one per node $v \in V$ (taking values in a finite set $E_v$), and a set of \emph{factors} $f_a: E_a \to \mathbb{R}$, one for each hyperedge $a \in \mathcal{H}E$. A joint distribution $P \in \mathbb{P}(E_V)$ associated with a factor graph satisfies the following decomposition: for any $x \in E_V$, as $P(x) = \prod_{a \in \mathcal{H}E} f_a(x_a)$.

\subsection{Generalized Bethe Free Energy: from Factor Graphs to presheaves over a Finite Poset with Values in Finite Sets}\label{presheaves}
\textbf{Bethe Free Energy for graphical models}

An undirected graphical model with respect to an acyclic graph, whose joint distribution $P\in \mathbb{P}_{>0}(E_V)$ is strictly positive (i.e., $P(x) > 0$ for $x \in E_v$), also factors according to its marginal distributions on the edges and vertices \cite{Speed, yedidia2005}. More precisely, consider an acyclic undirected graph $G = (V, E)$, i.e., there are no paths of consecutive edges starting at a vertex and ending at the same vertex. Let $P_a \in \mathbb{P}_{>0}(E_a)$ denote the marginal distribution of $X_a$; it is defined as $P_a(z_a) = \sum_{y_{\overline{a}} \in E_{\overline{a}}} P(z_a, y_{\overline{a}})$, where $z_a \in E_a$, and $\overline{a} = V \setminus a$ is the complement of $a$. Let $d(v)$ denote the degree of node $v \in V$, i.e., the cardinality of $\partial v$. Then, for any $x \in E_V$, 

\begin{equation}\label{graphical-model:factorization}
P(x) = \frac{\prod_{e \in E} P_e(x_e)}{\prod_{v \in V} P_v^{d(v) - 1}}.
\end{equation}

Recall that the entropy of a probability distribution $P\in \mathbb{P}(E_V)$ is defined as $S(P)=-\sum_{x\in E_V} P(x)\ln P(x)$. The factorization of Equation \ref{graphical-model:factorization} allows for rewriting the entropy of the joint distribution $S(P) = -\sum_{x \in E_V} P(x) \ln P(x)$ as

\begin{equation}\label{decomposition-entropy}
S(P) = \sum_{e \in E} S(P_e) - \sum_{v \in V} (d(v) - 1) S(P_v).
\end{equation}

 The previous decomposition of entropy for undirected graphical models on acyclic graphs (Equation \ref{decomposition-entropy}) is what motivates variational inference for undirected graphical models. Recall that in variational inference, one observes that for a joint distribution $P_{X,Y}$ over two variables $X \in E_X$ (the hidden variable) and $Y \in E_Y$ (the observed variable), the conditional probability $P_{X|Y}(x|y) = \frac{P_{X,Y}(x,y)}{P_Y(y)}$  given the observation $Y=y$ is the distribution $Q^{*} \in \mathbb{P}(E_X)$ that minimizes the free energy defined as $F(Q) = -\sum_{x \in E_X} Q(x) \ln P_{X,Y}(x,y) - S(Q)$.  In particular, the log-marginal distribution on the observed variables $-\ln P_Y(y)= -\ln \sum_{x\in E_X} P_{X,Y}(x,y)$ corresponds to the optimal value of $F(Q)$:

\begin{equation}
-\ln P_Y(y)= \inf_{Q\in \mathbb{P}(E_X)} F(Q)
\end{equation}
The Bethe free energy exploits the re-expression of entropy, Equation \ref{decomposition-entropy}, to express the free energy solely in terms of local marginals $Q_e \in \mathbb{P}_{>0}(E_e)$ for $e \in E$ and $Q_v \in \mathbb{P}_{>0}(E_v)$ for $v \in V$; for an undirected graphical model $P\in \mathbb{P}_{>0}(E_V)$ with respect to an acyclic graph $G=(V,E)$, the \emph{Bethe Free Energy}, denoted $F_{\text{Bethe}}$, is defined on the local marginals $(Q_e, Q_v; v\in V, e\in E)$ by the following expression:
\begin{equation}
F_{\text{Bethe}}\left(Q_a; a\in \mathcal{A}(G)\right) = \sum_{a \in \mathcal{A}(G)} c(a) \left[S(Q_a) - \mathbb{E}_{Q_a}[H_a]\right]
\end{equation}

with the following correspondence: $c(v) = -(d(v) - 1)$ for nodes, $c(e) = 1$ for edges, and $H_e = -\ln P_e$, $H_v = 0$. The expression of the Bethe free energy is also well-defined for undirected graphical models with graphs containing cycles; in this case, $F_{\text{Bethe}}(Q_a, a \in \mathcal{A}(G))$ differs from $F(Q)$ but is expected to be an approximation of $F(Q)$. Inference on undirected graphical models involves minimizing $F_{\text{Bethe}}\left(Q_a; a \in \mathcal{A}(G)\right)$, with $Q_a\in \mathbb{R}^{E_a}$ any collection of functions, under the constraints that for any $a \in \mathcal{A}(G)$, $Q_a$ is a stritcly probability distribution over $E_a$, i.e. $\sum_{x_a \in E_a} Q_a(x_a) = 1$ and for $x_a\in E_a$ $Q_a(x_a)>0$ , and the constraints that for any vertex $v \in V$, any edge $e = \{v, v_1\}$ such that $v_1 \in V$, and any $x_v \in E_v$, $Q_v(x_v) = \sum_{y \in E_{v_1}} Q_e(x_v, y)$. The second condition accounts for the fact that if there is a joint distribution $P \in \mathbb{P}(E_V)$ such that $Q_e$ is the marginal distribution of $P$ over $X_e$, i.e., $Q_e = P_e$, and similarly for vertices $Q_v = P_v$, then when $v \in e$, $Q_e$ must be related to $Q_v$ by marginalization.

More interestingly, following \cite{Yedidia}\cite{yedidia2005}\cite{PeltrePhD} such inference and the Bethe Free Energy extend to factor graphs and any collection of subsets of an index set of variables $I$. Their standard formulation is the basis of an extension of variational inference to presheaves from a poset to finite sets \cite{sergeantperthuis2022regionalized}\cite{sergeantperthuis:hal-04527780} which offers a framework general enough for us to then apply topological transformations to graphical models and trace their impact on the critical points of the associated Bethe Free Energies. In the next subsection we present both extensions; they rely on the inclusion-exclusion principle \cite{Rota}, which introduces weights $c(a)$ for subsets $a$ of an index set $I$. These weights $c(a)$ account for counting the same set multiple times due to overlaps and have a natural formulation for functions over a partially ordered set taking values in the set of real numbers, $\mathbb{R}$.

\textbf{Bethe Free Energy for factor graphs and extensions}

Let $\mathcal{A}$ be a finite partially ordered set. The ‘zeta-operator' of a poset $\mathcal{A}$, denoted $\zeta$, is the operator from $\bigoplus_{a \in \mathcal{A}} \mathbb{R} \to \bigoplus_{a \in \mathcal{A}} \mathbb{R}$, where $\bigoplus_{a \in \mathcal{A}} \mathbb{R}$ is identified with the set of functions from $\mathcal{A}$ taking values in $\mathbb{R}$, defined as, for any $\lambda \in \bigoplus_{a \in \mathcal{A}} \mathbb{R}$ and any $a \in \mathcal{A}$, $\zeta(\lambda)(a) = \sum_{b \leq a} \lambda_b$. It is a central result in combinatorics that the zeta-operator of $\A$ is invertible (Proposition 2 \cite{Rota}). We will call its inverse the M\"obius inversion of $\A$, and we will denote it as $\mu$. There is a collection $(\mu(a,b)\in \mathbb{Z}; b,a\in \A \text{ s.t. } b\leq a)$ such that, for any $\lambda\in \bigoplus_{a\in \A}\R$ and $a\in \A$, $\mu(\lambda)(a)=\underset{b\leq a}{\sum} \mu(a,b) \lambda_b$. The coefficients $c(a)$, which appear in the inclusion-exclusion principle, are defined as $c(a) = \sum_{b \geq a} \mu(b,a)$. The reader can refer to Appendix \ref{mobius-inversion-functor} for a detailed presentation of Möbius inversion and the inclusion-exclusion principle. Let $I$ be an index set of random variables, i.e., for each $i \in I$, there is a variable $X_i \in E_i$ taking values in a finite set $E_i$. Let $\mathcal{A}$ be a collection of subsets of $I$, i.e., $\mathcal{A} \subseteq \mathcal{P}(I)$ is a subset of $\mathcal{P}(I)$, the power set of $I$; in particular it is a poset for the inclusion relation, $\subseteq$, on subsets of $I$. Consider, for each $a \in \mathcal{A}$, a set of factors, $f_a: E_a \to \mathbb{R}$. Then the joint distribution $P \in \mathbb{P}(E_I)$ factors according to this collection of factors when, for any $x \in E_I$, $P(x) = \prod_{a \in \mathcal{A}} f_a(x_a)$. The Generalized Bethe free energy \cite{yedidia2005, PeltrePhD} associated with $\mathcal{A}$ and the collection $(f_a, a \in \mathcal{A})$ is defined as, for any collection of probability distributions $Q=(Q_a \in \mathbb{P}(E_a), a \in \mathcal{A})$,

\begin{equation}\label{bethe-free-energy}
F_{\text{Bethe}}(Q) = \sum_{a \in \mathcal{A}} c(a) \left( \mathbb{E}_{Q_a}[H_a] - S(Q_a) \right)
\end{equation}
with, for any $a \in \mathcal{A}$, $c(a) = \sum_{b \geq a} \mu(b,a)$, and for $x=(x_i,i\in a)$ an element of $E_a$, $H_a(x_a) = \sum_{b \subseteq a} -\ln f_b(x_b)$, where $x_b = (x_i, i \in b)$. We recover the case of graphical models when $\A=\A(G)$ for a graph $G$ and of factor graphs when $\A=\A(\mathcal{H})$ for a hypergraph $\mathcal{H}$. The collection of $Q_a, a \in \mathcal{A}$ that we are interested in are the ones compatible under marginalization, i.e., when $b \subseteq a$, we require that $Q_b$ is the marginal of $Q_a$ when summing out the variables $X_{a \setminus b}$. This constraint appears naturally when introducing presheaves that we shall define now. 

\begin{defn}[Graphical presheaves]
Let $I$ be a finite set, and let $\mathcal{A} \subseteq \mathcal{P}(I)$ be a sub-poset of the powerset of $I$. Let $E_i$, $i \in I$, be finite sets. For $a \in \mathcal{A}$, define $E_a := \prod_{i \in a} E_i$. Let $F_a := E_a$, and for $b \subseteq a$, let $F^a_b: E_a \to E_b$ be the projection map from $\prod_{i \in a} E_i$ to $\prod_{i \in b} E_i$. The functor $F$ is called a graphical presheaf from $\mathcal{A}$ to finite sets.

The probabilistic sections of a graphical presheaf are the collection of probability distributions $(Q_a \in \mathbb{P}(E_a), a \in \mathcal{A})$ such that, for any $b \subseteq a$,  
\[
\forall x \in E_b, \quad Q_b(x) = \sum_{y : F^a_b(y) = x} Q_a(y).
\]

\end{defn}

In the rest of the text, we shall rewrite the previous sum $Q_b(x) = \sum_{y \in F_a} 1[F^a_b(y) = x] Q_a(y)$, where $1[x \in A]$ is the indicator function which equals $0$ if $x \not\in A$ and $1$ if $x \in A$.

The probabilistic sections of a graphical presheaf $F$ are the sections, in the classical sense, of a presheaf that we will denote as $\mathbb{P}F$ and define just after. See Appendix \ref{sheaves} for a review of the classical definitions and constructions related to presheaves. Sections of presheaf $R$ correspond to elements $v_a \in R(a), \, a \in \mathcal{A}$ that satisfy the following compatibility property: whenever $b \leq a$, then $R^a_b(v_a) = v_b$. The set of sections of a presheaf $R$ is denoted by $\lim R$. The presheaf $\mathbb{P}F$ associates to each element $a \in \mathcal{A}$ the space of probability distributions $\mathbb{P}(E_a)$, viewed as a finite-dimensional vector space, where $E_a$ is the corresponding set. For inclusions $b \subseteq a$, it associates the marginalization map $\mathbb{P}F^a_b: \mathbb{P}(E_a) \to \mathbb{P}(E_b)$, which sends $P \in \mathbb{P}(E_a)$ to the pushforward measure defined as follows: for any $x \in E_b$, $\mathbb{P}F^a_b(P)(x) = \sum_{y \in E_a}1[ F^a_b(y) = x] P(y)$.

The optimization problem for variational inference graphical presheaves $F$, as introduced in \cite{yedidia2005, PeltrePhD}, involves minimizing the Bethe Free Energy $F_{\text{Bethe}}(Q_a; a \in \mathcal{A})$ under the convex constraint that the $(Q_a; a \in \mathcal{A})$ are probabilistic sections of $\mathbb{P}F$, i.e., $(Q_a; a \in \mathcal{A}) \in \lim \mathbb{P}F$. 

Graphical presheaves have the drawback of being very restrictive: the space $E_a$ must be a product of atomic sets $E_i$ for $i \in I$, with $i \in a$, corresponding to an `atomic' variable $X_i \in E_i$. Furthermore, the maps that relate values in $E_a$ to $E_b$, whenever \( b \subseteq a \), this must correspond to ignoring the variables \( X_{a \setminus b} \) following the projections from \( E_a = \prod_{i \in a} E_i \) to \( E_b = \prod_{i \in b} E_i \). These variables $X_a$ for $a \in \mathcal{A}$ correspond to the hidden variables that build the prior governing the dependencies among observed variables. It seems natural to rather associate with each $a \in \A$ a variable $Z_a \in F_a$ taking values in any finite set $F_a$. Furthermore, in this more general setting, one would want any possible relationship between $Z_a$ and $Z_b$, i.e., any function from $F_a$ to $F_b$ whenever $b\subseteq a$. In \cite{sergeantperthuis2022regionalized,sergeant-perthuis2024inference}, a generalization of the Bethe free energy and of the Belief Propagation algorithm are proposed for arbitrary variables $Z_a \in F_a$ indexed over any partially ordered set $(\mathcal{A}, \leq)$. The relationships between $Z_a$ and $Z_b$ when $b \leq a$ are given by stochastic maps, i.e., stochastic matrices whose rows consist of positive real values that sum to one. For the development of our paper, we only need the more restrictive case where the maps are deterministic, i.e., we allow any function $F^a_b: F_a \to F_b$ to relate the variables $Z_a$ and $Z_b$ when $b \leq a$. For this reason, in what follows, our main object of interest are presheaves from a poset $\mathcal{A}$ taking values in finite sets and the associated presheaf of probability distributions and their extension to probability distribution and vector spaces. We will denote \( \mathbf{FinSet} \) the category whose objects are finite sets and whose morphisms are maps between these sets (see Appendix~\ref{sheaves}).

\begin{defn}[Presheaf in finite sets and extensions]
Let $F$ be a presheaf over a finite poset $(\mathcal{A}, \leq)$ taking values in finite sets, i.e. a function that assigns to each element $a \in \mathcal{A}$ a finite set $F_a$ and a collection of functions $F^a_b: F_a \to F_b$ for any pair $b, a \in \mathcal{A}$ such that $b \leq a$, satisfying the condition that whenever $c \leq b \leq a$, $F^b_c \circ F^a_b = F^a_c$. 

We call the probabilistic extension of $F$ the presheaf $\mathbb{P}F$, which assigns to each element $a \in \mathcal{A}$ the set of probability distributions $\mathbb{P}(F_a)$. To any pair of elements $b, a \in \mathcal{A}$ such that $b \leq a$, it assigns the map $\mathbb{P}F^a_b: \mathbb{P}(F_a) \to \mathbb{P}(F_b)$, which sends a probability distribution $P \in \mathbb{P}(F_a)$ to its pushforward measure. For any $x \in F_b$, this is defined as  
\[
\mathbb{P}F^a_b(P)(x) = \sum_{y \in F_a}1[F^a_b(y) = x] P(y).
\]

Since $F_a$ is a finite set, $\mathbb{P}(F_a) \subseteq \mathbb{R}^{F_a}$ is a simplex in the vector space of real-valued functions $\mathbb{R}^{F_a}$. Furthermore, $\mathbb{P}F^a_b$ is the restriction of a linear map from $\mathbb{R}^{F_a}$ to $\mathbb{R}^{F_b}$. We denote the associated linear map by $\widetilde{F}^{a}_b: \mathbb{R}^{F_a} \to \mathbb{R}^{F_b}$, $\widetilde{F}^{a}_b(\lambda)(x)= \sum_{y \in F_a}1[F^a_b(y) = x] \lambda(y)$. The collection of spaces $(\widetilde{F}_a:= \mathbb{R}^{F_a};a\in \A)$ and of linear maps ($\widetilde{F}^a_b: \mathbb{R}^{F_a}\to \mathbb{R}^{F_b}; a,b\in \A \text{ such that } b\leq a)$, defines a presheaf over $\A$ taking values in finite vector spaces; we shall not it as $\widetilde{F}$. 

The dual map \( \widetilde{F}^{a,\ast}_b: \widetilde{F}_b^\ast \to \widetilde{F}_a^\ast \) sends linear forms \( l: \widetilde{F}_b \to \mathbb{R} \) to  $\widetilde{F}^{a,\ast}_b(l) = l \circ \widetilde{F}^a_b$.  
The collection of vector spaces \( (\widetilde{F}_a^\ast; a \in \mathcal{A}) \) and linear maps \( \widetilde{F}^{a,\ast}_b \) for \( b,a \in \mathcal{A}, b \leq a \) defines a functor over \( \mathcal{A} \) taking values in vector spaces. We shall denote this functor by \( \widetilde{F}^\ast \). When a scalar product \( \langle \cdot, \cdot \rangle_a \) is specified on each vector space \( \widetilde{F}_a \), there is an isomorphism $f_{\langle \cdot, \cdot \rangle_a}: F_a^{\ast} \to F_a$ 
which allows us to identify \( \widetilde{F}^{a,\ast}_b \) with its adjoint \( \widetilde{F}^{a,\dagger}_b \), defined as  $\widetilde{F}^{a,\dagger}_b = f_{\langle \cdot, \cdot \rangle_a} \circ \widetilde{F}^{a,\ast}_b \circ f_{\langle \cdot, \cdot \rangle_a}^{-1}$.

\end{defn}

The Generalized Bethe Free Energy for presheaves $F$ from a finite poset taking values in finite sets, and a collection of Hamiltonians $H_a: F_a \to \mathbb{R}$, has the same expression as Equation \ref{bethe-free-energy}, where $c(a) = \sum_{b \geq a} \mu(b,a)$ are defined using the Möbius inversion of the poset $\mathcal{A}$, and $Q$ are elements of $\lim F$.

\subsection{The General Belief Propagation Algorithm for Factor Graphs, Graphical Presheaves, and Generalization to presheaves of Sets}\label{belief-propagation}

The General Belief Propagation algorithm is used to find the critical points of the Generalized Bethe Free Energy. A classical result states that the fixed points of this algorithm correspond to the critical points of that free energy, which we state in Proposition~\ref{fix-points} and reprove in Appendix~\ref{proof-old}. Let us now recall the expression of this algorithm following \cite{yedidia2005,PeltrePhD}.

Let $I$ be a finite set that serves as index for variables $(X_i\in E_i; i\in I)$, each of which takes values in a finite set $E_i$; let $\A \subseteq \mathcal{P}(I)$ be a collection of subsets of I and denote $F$ the associated graphical presheaf. Let $(H_a:E_a\to \mathbb{R}, a \in \A)$ be a collection of Hamiltonians. Let us denote the update rule of the General Belief Propagation algorithm as $\text{BP}$. $\text{BP}$ acts on messages that we will now define. In the classical presentation of the algorithm, there are two types of messages at each time $t \in \mathbb{N}^{\ast}$. For elements $a, b \in \mathcal{A}$ such that $b \subseteq a$, we have top-down messages $m_{a \to b} \in \mathbb{R}_{>0}^{E_b}$ and bottom-up messages $n_{b \to a} \in \mathbb{R}_{>0}^{E_b}$. 

These messages are related as follows:

\begin{equation}\label{GBP1}
\forall a,b\in \A, \text{s.t. }  b\subseteq a, \quad n_{b \to a}^t= \prod_{\substack{c: b\subseteq c\\ c\not \subseteq  a}}m_{c\to b}^t
\end{equation}

Beliefs, which are interpreted as probability distributions up to a multiplicative constant and correspond to candidate sections of $F$, are defined as follows:

\begin{equation}\label{GBP2}
\forall a \in \A,\forall x_a\in E_a \quad b_a^t(x_a)\propto e^{-H_a(x_a)}\prod_{\substack{b\in \A:\\ b\leq a}} n_{b\to a}^t(x_b)
\end{equation}

where $\propto$ stands for proportional to.
The multiplication of function $n_{b\to a}$ that have different domains is made possible because there is an the embedding of $\R^{E_b}$ into $\R^{E_a}$ implicitly implied in the last equation; indeed, for $x \in E_a$ and $f \in E_b$, $f: x \mapsto f(x_b)$ defines a function from $E_a$ to $\mathbb{R}$.

To clarify the presentation, we require that $b_a$ be a probability distribution and normalize it accordingly. The update rule is given by,

\begin{equation}\label{GBP3}
\forall a,b\in \A, \text{s.t.} b\leq a\quad m_{a\to b}^{t+1}(x_b)= m_{a\to b}^{t}(x_b)\frac{\sum_{y_a: F^a_b(y_a)=x_b} b_a^t(y_a)}{b_b^t(x_b)}
\end{equation}

One observes that in the previous Equation~\ref{GBP3}, any normalization of beliefs does not change the update rule.
 The update rule can be rewritten in a more condensed manner, updating only the top-down messages, for all $a,b\in \A$, such that $b\leq a$,

\begin{align}\label{BP-update-rule}
 m_{a\to b}^{t+1}(x_b)&=  m_{a\to b}^t(x_b)\frac{\sum_{\substack{y_a: F^a_b(y_a)=x_b}} e^{-H_a(y_a)} \prod_{\substack{c\in \A:\\ c\subseteq a}}\prod_{\substack{d: c\subseteq  d\\ d\not \subseteq a}}m_{d\to c}^t (x_c)}{e^{-H_b(x_b)}\prod_{\substack{c\in \A:\\ c\subseteq b}}\prod_{\substack{d: c\subseteq  d\\ d\not \subseteq b}}m_{d\to c}^t (x_c)}
\end{align}

We denote the collection $(m_{a \to b}; a,b\in \A, b\leq a)$ as $m$. We denote $\text{BP}: \bigoplus_{a,b: b \subseteq a} \mathbb{R}^{E_b} \to \bigoplus_{a,b: b \subseteq a} \mathbb{R}^{E_b}$ as the operator underlying the update rule of Equation \ref{BP-update-rule}, i.e., we define $\text{BP}(m^t) = m^{t+1}$.

Consider the collection $\left(C_{a\to b} m_{a\to b}; a,b\in \A, b\leq a\right)$, where $C_{a\to b}$ is a strictly positive constant, i.e., it does not depend on $x_b\in F(b)$. Then, there is a collection of constants $(C'_{a\to b} > 0; a,b\in \A, b\leq a)$ such that  
\[
\text{BP}(C_{a\to b} m_{a\to b}) = C'_{a\to b} \text{BP}(m_{a\to b}).
\]

Furthermore, the associated beliefs defined by Equation~\ref{GBP2} remain unchanged under multiplication of $m_{a\to b}$ by a constant $C_{a\to b}$ for all $a,b\in \A$ such that $b\leq a$.

Therefore, $\text{BP}$ is an algorithm that preserves the equivalence classes $\{C \cdot m\}$, i.e., it is defined by the relation $m \sim m'$ whenever there is a collection of scalars $(C_{a\to b} \neq 0; a,b\in \A, b\leq a)$ such that  
\[
m_{a\to b} = C_{a\to b} m'_{a\to b} \quad \text{for any } a,b\in \A \text{ with } b\leq a.
\]
We shall denote the equivalence class of $m$ as $[m]$. The action of $\text{BP}$ on the equivalence classes of messages is denoted by $[\text{BP}]$ and defined as $[\text{BP}]([m]) = [\text{BP}(m)]$.

\begin{prop}[Yedidia, Freeman, Weiss, Peltre]\label{fix-points}
Let $I$ be an finite set, and $\A\subseteq \mathcal{P}(I)$ a collection of subsets of $I$; let $F$ be a graphical presheaf. Let $(m_{a\to b}\in \R_{>0}^{E_b},a,b\in \A \text{s.t. } b\subseteq a)$ be a fix point of the Generalized Belief Propagation up to a multiplicative constant, i.e. $[m]= [\text{BP}]([m])$. Let $(b_a,a\in \A)$ be the associated beliefs normalized so that each $b_a\in \mathbb{P}(E_a)$ (Equation \ref{GBP2}).
Then $(b_a,a\in \A)$ is a critical point of $F_{\text{Bethe}}$ under the constraint that $p\in \lim \mathbb{P}F$. Furthermore, any critical point of $F_{\text{Bethe}}$ in $\lim \mathbb{P}F$ is a belief associated to a fixed point of $\text{BP}$.
\end{prop}

\begin{proof}
See Theorem 5.15 in \cite{PeltrePhD} or Theorem 5 in \cite{yedidia2005}. This result is also a corollary of Theorem 2.2. \cite{sergeantperthuis2022regionalized}. In Appendix \ref{proof-old}, we reprove the result.

\end{proof}

An extension of the Belief Propagation algorithm to presheaves from a finite poset taking values in finite sets can be found in \cite{sergeantperthuis2022regionalized,sergeant-perthuis2024inference}. In their work, for each presheaf $F$, they propose a message-passing algorithm, that we will denote as $\text{MP}$, whose set of fixed points corresponds to the critical points of the Generalized Bethe Free Energy. When considering graphical presheaves, these message-passing algorithms slightly differ from the Belief Propagation algorithm but have the same fixed points, see Appendix~\ref{comparison-algorithm} for a description of the shared properties of the two algorithms. However, as we will show, $\text{MP}$ behaves well under natural transformations $\phi: F \to G$, whereas $\text{BP}$ does not. Let us now recall their message-passing algorithm, \text{MP}. To do so, following \cite{sergeantperthuis2022regionalized,sergeant-perthuis2024inference}, we need to introduce the elementary operators from which $\text{MP}$ is built.

For a functor $G$ from $\mathcal{A}$ to $\mathbb{R}$-vector spaces, we define $\mu_G$ as, for any $a\in \mathcal{A}$ and $v\in \bigoplus_{a\in \mathcal{A}} G_a$, $\mu_{G}(v)(a)=\sum_{b\leq a} \mu(a,b) G^b_a(v_b)$. Let $F$ be a presheaf from a finite poset $\A$ to finite sets $\textbf{FinSet}$. Let us call  $\text{FE}:\prod_{a\in\mathcal{A}}\R^{E_a}\to \prod_{a\in \mathcal{A}} \mathbb{R}$ the extension of $F_{\text{Bethe}}$ to real valued functions, i.e.  

\begin{equation}
\forall a \in \A, \text{FE}(h)(a)= \sum_{x_a\in F_a}h_a(x_a)H_a(x_a)+\sum_{x_a\in F_a} h_a(x_a)\ln h_a(x_a)
\end{equation}

Let us denote $b\leq a$ as $a\to b$. Let \(\delta_F: \bigoplus_{a \in \mathcal{A}} F_a \to \bigoplus_{a, b \in \mathcal{A} : b \leq a} F_b\) be defined as, for \((v_a \in F_a, a \in \mathcal{A})\):

\begin{equation}
\forall a, b \in \mathcal{A}, \text{s.t. } b \leq a, \quad \delta_F(v)(a \to b) = F^{a}_{b}(v_a) - v_b
\end{equation}

For a (covariant) functor \(G\), \(d_G: \bigoplus_{a, b \in \mathcal{A} : b \leq a} G(b) \to \bigoplus_{a \in \mathcal{A}} G_a\) is defined for \((v_{a \to b}, a, b \in \mathcal{A} \text{ such that } a \geq b)\) as

\begin{equation}
d_G(v)(a) = \sum_{b : b \leq a} G_{a}^{b}(v_{a \to b}) - \sum_{b : a \leq b} v_{b \to a}.
\end{equation}

Here, $G$ will be either $F^{\ast}$ or $F^{\dagger}$ when each space $F_a$ is equipped with a scalar product $\langle \cdot, \cdot \rangle_a$.

The $\zeta$ function of a functor $G$ plays an import role in the Belief propagation algorithm. $\zeta_G: \bigoplus_{a\in \mathcal{A}}G_a\to \bigoplus_{a\in \mathcal{A}}G_a $ is defined as, for $v\in \bigoplus_{a\in \mathcal{A}}G_a$, 

\begin{equation}
\zeta_{G}(v)(a)= \sum_{b\leq a} G^b_a(v_b)
\end{equation}

\begin{prop}\label{inverse-g-H}
Let $F$ be a presheaf from a finite to finite sets. Let $H_a\in \R^{F_a}, a\in \A$, be a collection of Hamiltonians. The differential $\text{d}\text{FE}: \bigoplus_{a\in \A}\widetilde{F}_a\to \bigoplus_{a\in \A}\widetilde{F}_a^\ast$ that sends $h$ to $\text{d}_h\text{FE}$ is invertible. We will denote its inverse by  $g_H: \bigoplus_{a \in \mathcal{A}} \widetilde{F}_a \to \bigoplus_{a \in \mathcal{A}} \widetilde{F}_a$.

\end{prop}
\begin{proof}

Recall that $\text{FE}_a(h_a) = \sum_{x_a} h_a(x_a) H_a(x_a) + \sum_{x_a} h_a(x_a) \ln h_a(x_a)$,
where $h \in \R^{E_a}$, and $\text{d}_x$ denotes the differential of the function at point $x$. Therefore,  
\begin{equation}\label{proof:differential}
\text{d}_h\text{FE}_a = \sum_{x_a} \text{d} h_a(x_a) (H_a(x_a) + \ln h_a(x_a) +1).
\end{equation}

with $\text{d} h_a(x_a) \in \bigoplus \widetilde{F}_a^\ast$ being the linear form acting as $(h_b(x_b); b\in \A, x_b\in F_b) \mapsto h_a(x_a)$. Consider the scalar products $\langle h_a, h^{'}_a\rangle_a=\sum_{x_a\in F_a}h_a(x_a) h^{'}_a(x_a)$ for each $a\in \A$; denote  $y_a\in \widetilde{F}_a$ the identification of $d_h\text{FE}_a\in F_a^{\ast}$ in $F_a$; then the reformulation of Equation \ref{proof:differential} is the following,

$$\forall x_a\in F_a, y_a(x_a)= H_a(x_a)+\ln h_a(x_a)+1 $$

The previous equation is equivalent to 
$$\forall x_a\in F_a, h_a(x_a)= e^{-H_a(x_a)+y_a(x_a)-1}$$

Therefore, we define  
\begin{equation}
\forall a\in \A, \forall x_a\in E_a, \quad g_{H_a,a}^{+}(l)(x_a) = e^{-H_a(x_a) + l_a(x_a) - 1}.
\end{equation}

The function $g_H = (g_{H_a,a}; a\in \A)$ is the inverse of $\text{d}\text{FE}$ when the image of $\text{d}\text{FE}$, i.e., $\bigoplus_{a\in \A} \widetilde{F}_a^\ast$, is identified with $\bigoplus_{a\in \A} \widetilde{F}_a$ through the scalar product  $\langle h, h_1 \rangle = \sum_{a\in \A} \langle h_a, h_{1,a} \rangle_a$.

\end{proof}

\begin{rem}
Any other choice of scalar product \( \langle\cdot, \cdot \rangle_1 \) on \( \widetilde{F}_a \) induces another map $g^{\langle\cdot, \cdot \rangle_1}_a: \widetilde{F}_a \to \widetilde{F}_a$ 
defined as $g_a \circ f_{\langle\cdot,\cdot\rangle} \circ f^{-1}_{\langle\cdot,\cdot\rangle_1}$, 
where \( \langle \cdot, \cdot \rangle \) is given by $\langle h, h_1 \rangle = \sum_{x_a \in E_a} h(x_a) h_1(x_a)$. 
The way the inverse  $g_H = \widetilde{F}_a^{\ast} \to \widetilde{F}_a$
of Proposition \ref{inverse-g-H} is represented as a function  $g_H^{\langle \cdot, \cdot \rangle}$ 
depends on the choice of scalar product on \( F_a \).  
For the scalar products $\langle h_a,h_{a}^{'}\rangle= \sum_{x_a} h_a(x_a)h_{a}^{'}(x_a)$,
we denote the identification $f^{\langle \cdot, \cdot \rangle}: \widetilde{F}_a^{\ast} \to \widetilde{F}_a$
as \( f_{+} \) and the associated inverse map as \( g_{H}^{+} \).
\end{rem}

In what follows, we denote $\bigoplus_{a,b \,:\, b \leq a} F_b$ as $F_{\rightarrow}$ to simplify the notations.

\begin{defn}[Message passing algorithms \cite{sergeantperthuis2022regionalized}]Choose for each \( a \in \mathcal{A} \) a scalar product \( \langle \cdot , \cdot \rangle_a \) on \( \widetilde{F}_a \); equip \( F_{\rightarrow} \) with the following scalar product: $\langle l, l_1 \rangle = \sum_{b,a \,:\, b \leq a} \langle l_{a\to b}, l_{1,a\to b} \rangle_b$.
The message-passing algorithm $\text{MP}_{F,H}$ for inference on a presheaf $F$ from a finite poset $\mathcal{A}$ to finite sets, with Hamiltonians $(H_a \in \mathbb{R}^{E_a}; a \in \mathcal{A})$, is defined as follows:

\begin{equation}
\forall l \in F_{\rightarrow}, \quad \text{MP}_{F,H}(l) = l+ \delta_F \circ g_H^{\langle \cdot,\cdot \rangle} \circ \zeta_{F^\dagger} \circ d_{F^\dagger}(l)
\end{equation}

The algorithm starts with a random initialization of messages $l_0 \in F_{\rightarrow}$, and $l^{t+1}$ is updated into the value $\text{MP}_{F,H}(l^{t})$. We will denote $\Delta MP_{F,H}$ the increment $\delta_F \circ g_H \circ \zeta_{F^\dagger} \circ d_{F^\dagger}$.

\end{defn}

\begin{prop}[Sergeant-Perthuis]\label{prop:fix-points}
Let $F$ be a presheaf from a poset $\mathcal{A}$ taking values in finite sets. Let $(H_a \in \mathbb{R}^{E_a}; a \in \mathcal{A})$ be a collection of Hamiltonians. Let $(l_{a \to b} \in \widetilde{F}_b^\ast; a,b\in \A: b \leq a)$ be a fixed point of $\text{MP}_{F,H}$, i.e., $\text{MP}_{F,H}(l) = l$. Let $b_a \in \mathbb{P}(F_a)$ be the unique probability distribution such that

$$b_a \propto g_H \circ \zeta_{F^\ast} \circ d_F(l).$$

Then $(b_a,a\in \A)$ is a critical point of $F_{\text{Bethe}}$ under the constraint that $p\in \lim \mathbb{P}F$. Furthermore, any critical point of $F_{\text{Bethe}}$ in $\lim \mathbb{P}F$ is a belief associated to a fixed point of $\text{MP}_{F,H}$. 

\end{prop}
\begin{proof}
Theorem 2.2. \cite{sergeantperthuis2022regionalized}. In Appendix \ref{proof-old}, we reprove the result.    
\end{proof}

\begin{rem}
Regardless of the choice of a scalar product, the fixed points of the associated message passing algorithm (Proposition \ref{prop:fix-points}) correspond to the same critical points. Indeed, the critical points are defined without the need to specify any scalar product. Furthermore, as explained in Appendix \ref{comparison-algorithm}, the fixed points of $\text{BP}$ and $\text{MP}$ are in correspondence; both algorithms (Propositions \ref{fix-points} and \ref{prop:fix-points}) allow us to characterize the critical points of the Bethe free energy.

\end{rem}

\section{Main results}\label{Main-results}

\subsection{Functoriality of $\Delta\text{MP}$}\label{functoriality}

Let $\A$ be a finite poset, and let $F, G$ be two presheaves taking values in finite sets. A natural transformation $\phi: F \to G$ is a collection of maps $\phi_a: F_a \to G_a$ compatible with the maps induced by the relation $b \leq a$, i.e., $G^a_b \circ \phi_a = \phi_b \circ F^a_b$ whenever $b \leq a$ (see Appendix \ref{sheaves}). A natural transformation from $F$ to $G$ induces a natural transformation between the presheaves $\mathbb{P}F$ to $\mathbb{P}G$ and $\widetilde{F}$ to $\widetilde{G}$ defined respectively as, for any $P\in \mathbb{P}(F_a)$ and $x\in G_a$, 
$\mathbb{P}\phi (P)(x)=\sum_{y\in F_a} 1[\phi_a(y)=x] P(y)$ and for any $h\in \mathbb{R}^{F_a}$, $\widetilde{\phi}(h)(x)=  \sum_{y\in F_a} 1[\phi_a(y)=x] h(y)$. Similarly, there is a natural transformation of functors \( \phi^{\ast}: \widetilde{G}^{\ast} \to \widetilde{F}^{\ast} \), i.e., a collection \( \phi^{\ast}_a: \widetilde{G}^{\ast}_a \to \widetilde{F}^{\ast}_a \). It is defined as follows: for any \( a \in \mathcal{A} \) and \( l \in \widetilde{G}_a^{\ast} \), $\phi^{\ast}_a(l) = l \circ \phi_a$. Similarly, when the spaces are equipped with a scalar product, there is a natural transformation $\phi^{\dagger}:\widetilde{G}^\dagger\to \widetilde{F}^\dagger$.

$\widetilde{\phi},\phi^{\ast},\phi^{\dagger}$ are extended into maps on the sources and targets of the operators $\delta_F$, $g_H$, $\zeta_{F^\ast}$, and $d_F$. $\phi$ extend to 
$
\phi^\oplus: \bigoplus_{a \in \mathcal{A}} \widetilde{F}_a \to \bigoplus_{a \in \mathcal{A}} \widetilde{G}_a$
where for $v\in \bigoplus_{a \in \mathcal{A}} \widetilde{F}_a$, $\phi^\oplus(v)_a= \widetilde{\phi}(v_a)$; we denote $\phi^\oplus$ simply as $\phi$. Similarly $\phi^{\oplus,*}: \bigoplus_{a \in \mathcal{A}} \widetilde{G}_a^{*} \to \bigoplus_{a \in \mathcal{A}} \widetilde{F}_a^{*}$, where $\phi^{*}(l)_a = l_a \circ \phi_a$; for the choice of scalar product on $\widetilde{F}_a$, $\phi^{\oplus,\dagger}=\bigoplus_{a\in \A} \phi_a^{\dagger}$. Similarly, \( \phi \) extends to \( \phi_1: F_{\rightarrow} \to G_{\rightarrow} \), where  $\phi(v)_{a\to b} = \widetilde{\phi}(v_{a\to b})$  
for \( v \in F_{\rightarrow} \). The map \( \phi_1^{\ast}: \widetilde{G}^\ast_{\rightarrow} \to \widetilde{F}^{\ast}_{\rightarrow} \) is defined as $\bigoplus_{a,b \,:\, b\leq a} \phi_b^\ast$ and $\phi_1^\dagger=\bigoplus_{a,b \,:\, b\leq a} \phi^\dagger_b$.\\

\begin{thm}\label{main:thm1}
Let \( F, G \) be two presheaves from a finite poset \( \mathcal{A} \) taking values in finite sets \textbf{FinSet}; let \( \phi: F \to G \) be a natural transformation. Then, for any collection of Hamiltonians \( (H_a: F_a \to \mathbb{R}; a \in \mathcal{A}) \), and any choice of scalar product on \( \widetilde{F}_a \) and \( \widetilde{G}_a \) with $a\in \A$,

\begin{equation}
\Delta MP_{G,\widetilde{H}}= \phi_1\circ \Delta MP_{F,H}\circ \phi_1^{\dagger}
\end{equation}

where for any $a\in \A$, 

\begin{equation}
\forall y\in G_a, \quad \widetilde{H}_a(x)=  -\ln \sum_{x\in F_a: \phi_a(x)= y}e^{-H_a(x)}
\end{equation}

with the convention that if \( x \not\in \operatorname{im} \phi_a \), then \( \widetilde{H}_a = +\infty \).

\end{thm}

\begin{proof}

In Appendix \ref{main-proofs}, we give a detailed proof of Theorem \ref{main:thm1}. It relies on showing the following properties: 
\begin{equation}
\begin{aligned}
\phi^{\ast} \circ d_G &= d_F \circ \phi_1^{\ast} \\
\phi^{\ast} \circ \zeta_{G^\ast} &= \zeta_{F^\ast} \circ \phi^{\ast}\\
\phi_1\circ \delta_F&= \delta_G\circ \phi
\end{aligned}
\end{equation}

One shows that for the choice of scalar product $\langle l, l_1\rangle_a= \sum_{x_a\in F_a} l(x_a)l_1(x_a)$, then,

$$\phi\circ g^{+}_{H}\circ \phi^{\dagger}= g^{+}_{\widetilde{H}}$$

Recall that \( f_{+} \) is the map that identifies \( F_a^\ast \) with \( F_a \), and \( g^{+}_{H} \) is the map associated with \( g_H \) by such an identification.

And one observes that for any other scalar product \( \langle \cdot,\cdot \rangle \),  
\[
\phi\circ g_H^{\langle \cdot,\cdot \rangle}\circ \phi^{\dagger}_{\langle \cdot,\cdot \rangle}= \phi\circ g^{+}_H\circ f_{+}\circ f_{\langle \cdot,\cdot\rangle}^{-1}\circ \phi^{\dagger}_{\langle\cdot,\cdot\rangle}
\]
as  
\[
g_H^{\langle\cdot,\cdot\rangle}= g_H^{+}\circ f^{+}\circ f_{\langle \cdot,\cdot\rangle}^{-1}.
\]
Then,  
\[
g_H^{\langle \cdot,\cdot \rangle}\circ \phi^{\dagger}_{\langle \cdot,\cdot \rangle}=  g_H^{+}\circ f_{+}\circ \phi^\ast\circ f^{-1}_{\langle \cdot,\cdot\rangle}= g_{H}^{+}\circ f_{+}\circ \phi^\ast \circ f_{+}^{-1}\circ f_{+}\circ f^{-1}_{\langle\cdot, \cdot\rangle}.
\]

As  
\[
\phi\circ g_H^{+}\circ f_{+}\circ \phi_\ast\circ f^{-1}_{+}= \phi \circ g_H^{+}\circ \phi^{\dagger}= g^{+}_{\widetilde{H}},
\]  
then,  
\[
\phi\circ g_H^{\langle \cdot,\cdot \rangle}\circ \phi^{\dagger}_{\langle \cdot,\cdot \rangle}= \phi\circ g^{+}_{\widetilde{H}}\circ f_{+}\circ f^{-1}_{\langle\cdot,\cdot\rangle}.
\]
Therefore,  
\[
\phi\circ g_H^{\langle \cdot,\cdot \rangle}\circ \phi^{\dagger}_{\langle \cdot,\cdot \rangle} = g_{\widetilde{H}}^{\langle \cdot,\cdot \rangle}.
\]

\end{proof}

\begin{rem}
By convention, we allow \( H_a(x_a) \) to be equal to \(+\infty \) for some \( a \in \mathcal{A} \) and \( x_a \in F_a \).  
If, furthermore, the collection of subsets \( S_a = \{x_a \mid H_a(x_a) = +\infty\} \) is such that $\overline{S_a}= F_a\setminus S_a$ defines a subobject of \( F \) that we will denote as $\overline{S}$, i.e. if $F^a_b(\overline{S_a}) \subseteq \overline{S_b}$, then for any $l\in F_{\rightarrow}$,
\[
\forall a, b, \, b \leq a, \quad \forall x_b \in S_b, \quad MP(l)_{a\to b}(x_b) = 0 
\]  

The proof of these statements can be written straightforwardly or can be seen as a consequence of Theorem~\ref{main:thm1} for \( \phi: \overline{S} \hookrightarrow F \), the inclusion of the subobject \( (\overline{S_a}; a \in \mathcal{A}) \) into \( F \). In particular, for any natural transformation \( \phi: F \to G \), the image \( \operatorname{im} \phi = (\operatorname{im} \phi_a; a \in \mathcal{A}) \) is a subobject of \( G \). To avoid having to handle points \( a \in \mathcal{A} \) and \( x_a \not\in \operatorname{im} \phi_a \) for which \( \tilde{H}_a(x_a) = +\infty \), one can restrict \( \text{MP} \) to messages \( l\in G_{\rightarrow} \) that satisfy \( l_{a \to b}(x_b) = 0 \) for all \( x_b \not\in \operatorname{im} \phi_b \).

Similarly, one can choose as a convention that \( BP(m)_{a\to b}(x_b) = -\infty \) for $x_b\in S_b$; more precisely, in both case for $\text{MP}$ and $\text{BP}$ both conventions imply that the associated beliefs are \( 0 \):  
\[
b_a(x_a) = g_H^{\langle\cdot,\cdot\rangle} \circ \zeta_{F^{\dagger}} \circ d_{F^{\dagger}}(\ln m)(a, x_a) = 0,
\]
Then for \( x_b \in S_b \), we have  
\[
F^a_b (b_a)(x_b) = 0.
\]

\end{rem}

\begin{prop}[The category \textbf{MPA}]
    Let us now define the category of message passing algorithms, which we will denote as $\textbf{MPA}$. Its objects are maps \( D_0:H \times S \to H \), where \( H \) is a Hilbert space and \( S \) is a topological space; those maps correspond to dynamics from a Hilbert space \( H \) to itself, indexed by a parameter space \( S \). 

Consider a second object, \( D_1:H_1 \times S_1 \to H_1 \). A morphism between two such objects, \( \phi: D_0 \to D_1 \), is a pair \( (\phi_0, \phi_1) \), where \( \phi_0: H \to H_1 \) is a bounded linear map and \( \phi_1: S \to S_1 \) is a continuous map, such that

\[
\forall x \in H_1, \forall s \in S, \quad D_1(x, \phi_1(s)) = \phi_0 \circ D_0(\phi_0^{\dagger}(x), s).
\]

The composition \( \phi \circ \psi \) of two morphisms \( \phi, \psi \) is defined as the composition of their maps \( \phi_0 \circ \psi_0 \) and \( \phi_1 \circ \psi_1 \).
\end{prop}
\begin{proof}

 \( \phi_0 \circ \psi_0 \) and \( \phi_1 \circ \psi_1 \) are respectively bounded linear and convex; the composition is indeed a morphism, $\phi\circ \psi: D_0\to D_2$,  from $D_0$ to $D_2$ as:

\[
\forall x \in H_2, \forall s \in S, \quad D_2(x, (\phi_1 \circ \psi_1)(s)) = \phi_0 \circ D_1(\phi_0^{\dagger}(x), \phi_1(s)) = \phi_0 \circ \psi_0 \circ D_0(\psi_0^{\dagger} \circ \phi_0^{\dagger}(x), s).
\]
and so,
\[
\forall x \in H_2, \forall s \in S, \quad D_2(x, (\phi_1 \circ \psi_1)(s))= \phi_0 \circ \psi_0 \circ D_0((\phi_0 \circ \psi_0)^\dagger x, s)
\]

\end{proof}

Let us denote the category of presheaves from a poset \( \mathcal{A} \) to finite sets as \( \hat{\mathcal{A}}_f \). The space \( \prod_{a \in \mathcal{A}} \mathbb{R}^{F_a} \) is a finite-dimensional vector space and therefore a topological space for any norm defined on it. Consider, for example, the \( \|\cdot\|_{\infty} \) norm, defined as  $\| H \|_{\infty} = \sup_{a \in \mathcal{A}} \sup_{x_a \in F_a} |H_a(x_a)|$.

\begin{thm}\label{fonctoriality-increment}
Let $\A$ be a finite poset; let $F:\A\to \textbf{FinSet}$ be a presheaf from $\A$ to finite sets. Recall that $F_{\rightarrow}=\bigoplus_{a,b\in \A: b\leq a} F_a$; equip each $\tilde{F}_a=\R^{F_a}$ with any scalar product $\langle \cdot, \cdot \rangle_a$ which in turn equips $F_\rightarrow$ with a scalar product. The map $\Delta MP(F): F_{\rightarrow}\times \prod_{a\in \A}\R^{F_a}\to F_{\rightarrow} $ is  defined as, 

$$\forall l \in F_{\rightarrow}, H\in\prod_{a\in \A}\R^{F_a},\quad  \Delta MP(F)(l, H)=\Delta MP_H(l) $$

To a natural transformation \( \phi: F \to G \) between two presheaves of \( \hat{\mathcal{A}}_f \), \( \Delta MP(\phi) = (\phi, \psi) \) is a morphism between \( \Delta MP(F) \) and \( \Delta MP(G) \) in \( \text{MPA} \), where  

\[
\forall a \in \mathcal{A}, \, y_a \in G_a, \quad \psi(H)_a(y_a) = -\ln \sum_{x_a \in F_a \,:\, \phi_a(x_a) = y_a} e^{-H_a(x_a)}
\]

Furthermore, \( \Delta MP \) is a functor from $\hat{\A}_f\to \textbf{MPA}$.

\end{thm}

\begin{proof}
It is a direct consequence of Theorem \ref{main:thm1}. See Appendix \ref{main-proofs} for the full proof.
\end{proof}

\subsection{Functoriality of $\text{MP}$}

Let \( F \) be a presheaf of \( \hat{\mathcal{A}}_f \), and recall that for \( m \in F_{\rightarrow} \), \( \text{MP}_{F,H}(m) = m + \Delta \text{MP}_{F,H}(m) \). Let \( G \) be a second presheaf of \( \hat{\mathcal{A}}_f \), and let \( \phi: F \to G \) be a natural transformation. Then, by Theorem \ref{main:thm1}, \( \text{MP}_{G,\tilde{H}}(m) = m + \phi \circ \Delta\text{MP}_{F,H}(\phi^\dagger(m)) \); however, \( \phi \circ \phi^\dagger \) is not necessarily equal to the identity map, and therefore \( \text{MP}_{G,\tilde{H}} \) is, in general, different from \( \text{MP}_{G,\tilde{H}} \circ \phi^\dagger \). We will show in this section that for a good choice of scalar products on the $F_a$'s, with $a\in \A$, $\text{MP}_{F,H}$ is also a functor. 

\begin{lem}\label{lem:isometry}
Let $F, G$ be two presheaves from a poset $\mathcal{A}$ to a finite measurable set, let $\phi: F\to G$ be a natural transformation. Consider on $\widetilde{G}_a$ the standard scalar product $\langle h,j\rangle= \sum_{x\in G_a} h(x_a)j(x_a)$; for any $a \in \mathcal{A}$, consider the following scalar product on $\widetilde{F}_a$:
\[
\langle f, g \rangle_a^\phi = \sum_{x_a \in F_a} \vert \phi_a^{-1}(\{ \phi_a (x_a)\} ) \vert f(x_a) g(x_a).
\]

Then \( \phi_a: \widetilde{F}_a \to \widetilde{G}_a \) and its adjoint \( \phi^{\dagger} \) satisfy \( \phi \circ \phi^{\dagger}(f) = f \) for any \( f \in \prod_{a \in \mathcal{A}} \mathbb{R}^{\im \phi_a} \), i.e., \( \phi \circ \phi^{\dagger} \) is the identity map on \( \widetilde{\im \phi} \); furthermore, \( \phi^{\dagger} \) is an isometry. In particular, if \( \phi \) is a monomorphism, i.e. if for all $a\in \A$ $\phi_a$ is surjective, then \( \phi \circ \phi^\dagger = \text{id} \). 
 
\end{lem}

\begin{proof}
Let $f\in \widetilde{G}_a$ and $g\in \widetilde{F}_a$, then,

$$ \langle f, \phi(g) \rangle= \sum_{y_a\in G_a} f(x_a) \sum_{x_a:\phi_a(x_a)=y_a} g(x_a)$$
$$=\sum_{x_a\in F_a} f(x_a) g(\phi_a(x_a))\frac{\vert \phi_a^{-1}(\{ \phi_a (x_a)\} ) \vert}{\vert \phi_a^{-1}(\{ \phi_a (x_a)\} ) \vert}  $$

Therefore, $\phi^{\dagger}(g)(x_a)= \frac{1}{\vert \phi_a^{-1}(\{ \phi_a (x_a)\} ) \vert}g\circ \phi(x_a)$. Furthermore, when $\phi_a:F_a\to G_a$ is surjective,
\[
\forall y_a \in G_a, \quad (\phi \circ \phi^{\dagger} g)(y_a) = \sum_{x_a : \phi(x_a) = y_a} \frac{1}{\vert \phi_a^{-1}(\{ \phi_a (x_a)\} ) \vert} g \circ \phi(x_a).
\]
When $\phi_a^{-1}(\{y_a\})$ is empty, the previous expression is $0$. If it is not empty—which is the case when $\phi: F_a \to G_a$ is surjective—then

\[
\sum_{x_a : \phi(x_a) = y_a} \frac{1}{\vert \phi_a^{-1}(\{ \phi_a (x_a)\} ) \vert} g \circ \phi(x_a) = \sum_{x_a : \phi(x_a) = y_a} \frac{1}{\vert \phi_a^{-1}(\{ y_a\}) \vert} g \circ \phi (x_a) = g(y_a).
\]

In particular for $f,g\in \tilde{\im \phi}_a$, one has that,
\[
\langle \phi^{\dagger}(f), \phi^{\dagger}(g) \rangle^\phi_a = \langle \phi \circ \phi^{\dagger} f, g \rangle_a = \langle f, g \rangle_a,
\]
which shows that $\phi^{\dagger}$ is an isometry.

\end{proof}

\begin{thm}\label{label:thm-3}
Let \( \mathcal{A} \) be a finite poset; let \( F,G: \mathcal{A} \to \textbf{FinSet} \) be two presheaves from \( \mathcal{A} \) to finite sets. Let $\phi: F\to G$ be a monomorphisms, i.e. a natural transformation made of surjective maps. Recall that \( F_{\rightarrow} = \bigoplus_{a,b \in \mathcal{A}: b \leq a} F_a \); equip each \( \widetilde{F}_a = \mathbb{R}^{F_a} \) with any scalar product \( \langle \cdot, \cdot \rangle_a^\phi \), which in turn equips \( F_\rightarrow \) with a scalar product. Let for $l\in F_\rightarrow$ and $H\in \prod_{a\in \A}\R^{F_a}$, \( \text{MP}(F)(l, H) = l+\delta_{F} \circ g_{H}^{\langle \cdot, \cdot \rangle^\phi} \circ \zeta_{F^\dagger} \circ d_{F^\dagger}(l) \); let \( \text{MP}(\phi) = (\phi, \psi) \) with,

\[
\forall a \in \mathcal{A}, \, y_a \in G_a, \quad \psi(H)_a(y_a) = -\ln \sum_{x_a \in F_a : \phi_a(x_a) = y_a} e^{-H_a(x_a)}.
\]

Then \( \text{MP} \) is a functor from \( \hat{\mathcal{A}}_f \) to \( \textbf{MPA} \).
\end{thm}
\begin{proof}
Theorem \ref{label:thm-3} is a direct consequence of Lemma \ref{lem:isometry} and Theorem \ref{fonctoriality-increment} as for any $l\in G_\rightarrow$,

$$
\phi \circ \text{MP}_{F,H} \circ \phi^{\dagger}(l) = \phi \circ \phi^\dagger(l) + \phi \circ \text{MP}_{F,H} \circ \phi^{\dagger}(l)
= \text{MP}_{G,\widetilde{H}}(l)$$
For a detailed proof of the result, see Appendix \ref{main-proofs}.

\end{proof} 
\section{Discussion}

The Belief Propagation algorithms and associated message passing algorithms are in correspondence with fixed points of the Bethe Free Energy, which allow for approximate inference on factor graphs and generalization. Exploiting recent results relating such algorithms to presheaves from a poset to finite sets, we propose a first systematic study of the impact of topological transformations on these message passing algorithms. We show that natural transformations of those presheaves relate to linear maps on the associated message passing algorithms, complementing results on the impact of the topology of the factor graph on those algorithms \cite{welling2012choice, PeltrePhD}. Such novel results are a first step towards a unified study on the role the topology of factor graphs plays in approximate inference of factor graphs.

\section{Acknowledgment}
The authors would like to thank the PRMO Young FMJH Program for supporting this research.

\bibliographystyle{plain}
\bibliography{hal/hal-bib}

\begin{thebibliography}{10}

\bibitem{barbero2022sheaf}
Federico Barbero, Cristian Bodnar, Haitz~S{\'a}ez de~Oc{\'a}riz~Borde, Michael Bronstein, Petar Veli{\v{c}}kovi{\'c}, and Pietro Li{\`o}.
\newblock Sheaf neural networks with connection laplacians.
\newblock In {\em Topological, Algebraic and Geometric Learning Workshops 2022}, pages 28--36. PMLR, 2022.

\bibitem{Baudot2015}
Pierre Baudot and Daniel Bennequin.
\newblock The homological nature of entropy.
\newblock {\em Entropy}, 17(5):3253--3318, 2015.

\bibitem{bishop2006pattern}
Christopher~M. Bishop.
\newblock {\em Pattern Recognition and Machine Learning}.
\newblock Springer, 2006.

\bibitem{10.5555/2500991}
Alain Bretto.
\newblock {\em Hypergraph Theory: An Introduction}.
\newblock Springer Publishing Company, Incorporated, 2013.

\bibitem{Carlsson2006}
Erik Carlsson, Gunnar Carlsson, and Vin De~Silva.
\newblock An algebraic topological method for feature identification.
\newblock {\em International Journal of Computational Geometry \& Applications}, 16(04):291--314, 2006.

\bibitem{Curry2013}
Justin Curry.
\newblock {\em Sheaves, cosheaves and applications}.
\newblock PhD thesis, {T}he {U}niversity of {P}ennsylvania, 2013.
\newblock arXiv:1303.3255.

\bibitem{gallier2022homology}
Jean Gallier and Jocelyn Quaintance.
\newblock {\em Homology, Cohomology, and Sheaf Cohomology for Algebraic Topology, Algebraic Geometry, and Differential Geometry}.
\newblock World Scientific, 2022.

\bibitem{Ghrist2008}
Robert Ghrist.
\newblock Barcodes: the persistent topology of data.
\newblock {\em Bulletin of the American Mathematical Society}, 45(1):61--75, 2008.

\bibitem{gimpel2008statistical}
Kevin Gimpel and Daniel Rudoy.
\newblock {\em Statistical inference in graphical models}.
\newblock Massachusetts Institute of Technology, Lincoln Laboratory, 2008.

\bibitem{grim2024message}
Anna Grim.
\newblock {\em Message Passing Dynamics of Belief Propagation Algorithms}.
\newblock Doctoral dissertation, Brown University, May 2022.

\bibitem{hansen2020sheafneuralnetworks}
Jakob Hansen and Thomas Gebhart.
\newblock Sheaf neural networks, 2020.

\bibitem{8919796}
Jakob Hansen and Robert Ghrist.
\newblock Distributed optimization with sheaf homological constraints.
\newblock In {\em 2019 57th Annual Allerton Conference on Communication, Control, and Computing (Allerton)}, pages 565--571, 2019.

\bibitem{Kashiwara}
Masaki Kashiwara and Pierre Schapira.
\newblock Persistent homology and microlocal sheaf theory.
\newblock {\em Journal of Applied and Computational Topology}, 2018.

\bibitem{Lauritzen}
Steffen~L. Lauritzen.
\newblock {\em Graphical Models}.
\newblock Oxford Science Publications, 1996.

\bibitem{lauritzen2019lecture}
Steffen~L. Lauritzen.
\newblock Lectures on graphical models.
\newblock Lecture Notes, 2019.

\bibitem{mac2013categories}
Saunders Mac~Lane.
\newblock {\em Categories for the working mathematician}, volume~5.
\newblock Springer Science \& Business Media, 2013.

\bibitem{10.5555/1592967}
Marc Mezard and Andrea Montanari.
\newblock {\em Information, Physics, and Computation}.
\newblock Oxford University Press, Inc., USA, 2009.

\bibitem{Pearl1988}
Judea Pearl.
\newblock {\em Probabilistic Reasoning in Intelligent Systems: Networks of Plausible Inference}.
\newblock Morgan Kaufmann Publishers, 1988.

\bibitem{Peltre}
Olivier Peltre.
\newblock {Homology of Message-Passing Algorithms}.
\newblock \url{http://opeltre.github.io}, 2020.
\newblock Ph.D. thesis (preprint).

\bibitem{PeltrePhD}
Olivier Peltre.
\newblock Message passing algorithms and homology, 2020.
\newblock Ph.D. thesis, \href{https://opeltre.github.io/assets/bib/Peltre-Message_Passing_Algorithms_and_Homology.pdf}{Link}.

\bibitem{riehl2017category}
Emily Riehl.
\newblock {\em Category theory in context}.
\newblock Courier Dover Publications, 2017.

\bibitem{Rota}
Gian-Carlo Rota.
\newblock On the foundations of combinatorial theory {I}. {Theory} of {M}{\"o}bius functions.
\newblock {\em Probability theory and related fields}, 2(4):340--368, 1964.

\bibitem{sergeantperthuis:hal-04527780}
Gr{\'e}goire Sergeant-Perthuis and Nils Ruet.
\newblock {Inference on diagrams in the category of Markov kernels (Extended abstract)}.
\newblock In {\em {7th International Conference on Applied Category Theory (ACT 7)}}, Oxford (UK), United Kingdom, June 2024. {David Jaz Myers and Michael Johnso}.

\bibitem{sergeantperthuis2022regionalized}
Grégoire Sergeant-Perthuis.
\newblock Regionalized optimization, 2022.

\bibitem{sergeant-perthuis2024inference}
Grégoire Sergeant-Perthuis and Nils Ruet.
\newblock Inference on diagrams in the category of markov kernels (extended abstract).
\newblock 2024.
\newblock Presented at ACT 7.

\bibitem{Speed}
Terry~P. Speed.
\newblock A note on nearest-neighbour gibbs and markov probabilities.
\newblock {\em Sankhy\=a: The Indian Journal of Statistics, Series A}, 1979.

\bibitem{Vigneaux2020information}
Juan~Pablo Vigneaux.
\newblock Information structures and their cohomology.
\newblock {\em Theory and Applications of Categories}, 35(38):1476--1529, 2020.

\bibitem{wainwright2008graphical}
Martin~J Wainwright, Michael~I Jordan, et~al.
\newblock Graphical models, exponential families, and variational inference.
\newblock {\em Foundations and Trends{\textregistered} in Machine Learning}, 1(1--2):1--305, 2008.

\bibitem{welling2012choice}
Max Welling.
\newblock On the choice of regions for generalized belief propagation.
\newblock {\em arXiv preprint arXiv:1207.4158}, 2012.

\bibitem{10.5555/3008751.3008848}
Jonathan~S. Yedidia, William~T. Freeman, and Yair Weiss.
\newblock Generalized belief propagation.
\newblock In {\em Proceedings of the 13th International Conference on Neural Information Processing Systems}, NIPS'00, page 668–674, Cambridge, MA, USA, 2000. MIT Press.

\bibitem{Yedidia}
Jonathan~S Yedidia, William~T Freeman, and Yair Weiss.
\newblock Constructing free-energy approximations and generalized belief propagation algorithms.
\newblock {\em IEEE Transactions on information theory}, 51(7):2282--2312, 2005.

\bibitem{yedidia2005}
J.S. Yedidia, W.T. Freeman, and Y.~Weiss.
\newblock Constructing {Free}-{Energy} {Approximations} and {Generalized} {Belief} {Propagation} {Algorithms}.
\newblock {\em IEEE Transactions on Information Theory}, 51(7):2282--2312, July 2005.

\end{thebibliography}

\appendix

\section{Categories, posets, functors and natural transformations}\label{sheaves}

For an introduction to categories, see \cite{mac2013categories,riehl2017category}. A category $\textbf{C}$ is a collection of objects \( \mathcal{O}(\textbf{C}) \) and morphisms \( \phi: A \to B \) between two objects. Two morphisms \( \phi: A \to B \), \( \phi_1: B \to C \), whose target and domain are compatible, compose into a morphism \( \phi_1 \circ \phi: A \to C \). The composition is associative, \( \phi_2 \circ (\phi_1 \circ \phi) = (\phi_2 \circ \phi_1) \circ \phi \). We will denote categories in bold. A (covariant) functor $F:\textbf{C}\to \textbf{C}_1$ between two categroies $\textbf{C}$ and $\textbf{C}_1$ send objects $A$ of $\mathcal{O}(\textbf{C})$ to objects $F(A)$ of $\mathcal{O}(\textbf{C}_1)$; it sends a morphism $\phi: A\to B$ to a morphism $F(\phi): F(A)\to F(B)$ such that $F(\phi_1\circ \phi)= F(\phi_1)\circ F(\phi)$. A contravariant functor sends morphisms \( \phi \) to \( F(\phi): F(B) \to F(A) \) and satisfies \( F(\phi_1 \circ \phi) = F(\phi) \circ F(\phi_1) \).  
A covariant functor \( F \) from \( \textbf{C} \) to \( \textbf{C}_1 \) is denoted as \( F:\textbf{C} \to \textbf{C}_1 \), and a contravariant functor as \( F:\textbf{C}^{op} \to \textbf{C}_1 \).

The categories of sets and finite sets, denoted respectively by \( \textbf{Set} \) and \( \textbf{FinSet} \), have as objects sets and finite sets, respectively, and as morphisms between two sets \( X, Y \), functions from \( X \) to \( Y \). A contravariant functor from a category \( \textbf{C} \) to \( \textbf{Set} \) is called a presheaf. In particular, contravariant functors from \( \textbf{C} \) to \( \textbf{FinSet} \) are presheaves.
We will also encounter implicitly the categroy of finite vector spaces, which has as objects finite vector spaces and as morphims linear maps between two such vector space. 

A natural transformation between two covariant or contravariant functors \( F,G:\textbf{C} \to \textbf{C}_1 \) or \( F,G:\textbf{C}^{op} \to \textbf{C}_1 \)  
is a collection of morphisms \( \phi_A: F(A) \to G(A) \) indexed by objects \( A \) of \( \mathcal{O}(\textbf{C}) \)  
such that for any morphism \( \psi: A \to B \), we have $\phi_B \circ F(\psi) = G(\psi) \circ \phi_A$ when \( F \) and \( G \) are covariant functors, and $\phi_A \circ G(\psi) = F(\psi) \circ \phi_B$. A natural transformation \( \phi \) between two covariant or contravariant functors \( F,G:\textbf{C} \to \textbf{Set} \)  
is called a monomorphism if each \( \phi_A \) is injective, and an epimorphism if each \( \phi_A \) is surjective.

A partially ordered set (poset), denoted as \( \A \), is a set equipped with a binary relation \( \mathcal{R} \subseteq \A \times \A \) which satisfies:

\begin{enumerate}
    \item Reflexivity: \( \forall x \in \A, \quad x \mathcal{R} x \).
    \item Transitivity: \( \forall x,y,z \in \A \), if \( x \mathcal{R} y \) and \( y \mathcal{R} z \), then \( x \mathcal{R} z \).
    \item Antisymmetry: \( \forall x,y \in \A \), if \( x \mathcal{R} y \) and \( y \mathcal{R} x \), then \( x = y \).
\end{enumerate}

$\mathcal{R}$ is usually dentoted as $\leq$.
One associates to a poset a category \( \mathcal{C}(\A) \), whose objects \( \mathcal{O}(\A) \) are the elements of \( \A \), and there is a unique morphism \( \phi: b \to a \) between \( a, b \in \A \) when \( b \leq a \). We will denote the associated category \( \mathcal{C}(\A) \) simply as \( \A \). A contravariant functor \( F: \A \to \textbf{C} \) from the poset \( \A \) seen as a category to a category \( \textbf{C} \) is encoded into a collection of objects \( (F_a, a \in \A) \) and a collection of morphisms \( (F^a_b: F_a \to F_b; \, a,b \in \A, \, b \leq a) \).

The set of sections of a presheaf $F$, from a finite poset $\A$ to $\textbf{Set}$, also called limit and denoted $\lim F$, is the set of collections of elements $x = (x_a \in F_a \mid a\in \A)$ such that they are pairwise compatible: for any $a,b\in \A$ with $b\leq a$, $F^a_b(x_a) = x_b$.

The application of presheaves and sheaves in data science has gained more importance in recent years; they were introduced in the context of decentralized optimization, for deep learning architectures that enable heterogeneous descriptions of data (Sheaf neural networks) \cite{Curry2013,8919796,hansen2020sheafneuralnetworks,barbero2022sheaf}. These approaches involve introducing sheaves over cell complexes as a data structure \cite{Curry2013} and using a Laplacian as the foundational building block for neural network layers. A generalization of inference on graphical models, which involves presheaves from a poset to the category of Markov kernels, was introduced in \cite{sergeant-perthuis2024inference}, and for decentralized optimization in \cite{sergeantperthuis2022regionalized}. They also appear in topological data analysis \cite{Carlsson2006, Ghrist2008,Kashiwara} and for information theory \cite{Baudot2015,Vigneaux2020information}.

\section{Möbius inversion of a poset and of a functor}\label{mobius-inversion-functor}

The definition of the Bethe Free Energy relies on weights given by the inclusion-exclusion formula of a poset; we shall now recall its definition.

\begin{defn}[Zeta operator of a poset]
Let $\A$ be a finite poset. We call the ‘zeta-operator' of a poset $\A$, denoted $\zeta$, the operator from $ \bigoplus_{a\in \A}\R\to \bigoplus_{a\in \A}\R$ defined as, for any $\lambda\in \bigoplus_{a\in \A}\R$ and any $a\in \A$,
\begin{equation}
\zeta(\lambda)(a)=\underset{b\leq a}{\sum} \lambda_b
\end{equation}

\end{defn}

\begin{prop}[Reformulation of Proposition 2 \cite{Rota}, Rota' 64]\label{mobius-inversion}
Let $\A$ be a finite poset. The zeta-operator of $\A$ is invertible. We will call its inverse the M\"obius inversion of $\A$, denoted $\mu$. Furthermore, there is a collection $(\mu(a,b); b,a\in \A \text{ s.t. } b\leq a)$ such that, for any $\lambda\in \bigoplus_{a\in \A}\R$ and $a\in \A$,

\begin{equation}
\mu(\lambda)(a)=\underset{b\leq a}{\sum} \mu(a,b) \lambda_b
\end{equation}

\end{prop}

We call the coefficient $(\mu(a,b), b,a \text{ s.t. } b\leq a)$ the  M\"obius coefficients of $\A$. In particular Proposition \ref{mobius-inversion} implies that, for any $b,a\in \A$ such that $b\leq a$,

\begin{align}
\sum_{c:b\leq c\leq a} \mu(a,c) &= 1[b=a]\\
\sum_{c:b\leq c\leq a} \mu(c,b) &= 1[b=a]
\end{align}

For a collection of values $\lambda_a\in \mathbb{R}$, $a\in \mathcal{A}$, we call the following expression $\sum_{a\in \mathcal{A}}\sum_{b\leq a} \mu(a,b) \lambda_b$ the inclusion-exclusion formula over a poset $\mathcal{A}$; this formula corresponds to the value one would attribute to a maximal element, denoted $1$, added to $\mathcal{A}$:

\begin{equation}
\lambda_{1}:=\sum_{a\in \A}\sum_{b\leq a} \mu(a,b) \lambda_b
\end{equation}

which can be rewritten as $\lambda_{1}=\sum_{a\in \A}[\sum_{b\geq a} \mu(b,a)] \lambda_a$. We will denote $c(a)=\sum_{b\leq a} \mu(a,b)$ the weighted coefficients.

To find the classical inclusion-exclusion formula, consider $I$ a finite set; the poset $\mathcal{A}$ is $(\mathcal{P}(I),\supseteq)$ with the reversed order. The quantities $\lambda_a$, for $a \in \mathcal{A}$, are $\vert A_i \vert$ when $a=i \in I$, and represent the cardinality of the sets $A_i$, as well as the cardinality of all possible intersections $\vert \cap_{i\in a} A_i \vert$ when $a\subseteq I$. In this setting, the maximal element $1$ has a value $\lambda_1 = \vert \cup_{i\in I} A_i \vert$.\\

Let $G$ be a (covariant) functor from a finite poset $\A$ to finite vector spaces. Define,

$$\forall v\in \bigoplus G_a, \quad \zeta_{G}(v)_a=\sum_{b\leq a}G^{b}_a(v_b)$$

and,

$$\forall v\in \bigoplus G_a, \quad \mu_{G}(v)_a=\sum_{b\leq a}\mu(a,b)G^{b}_a(v_b).$$

Then, 

$$\eta_G\circ \mu_G(v)_a=\sum_{c,b: c\leq b\leq a} \mu(b,c) G^c_a(v_c)= v_a,$$

and similarly, $\mu_G\circ \eta_G=\operatorname{id}$.

\section{Proof of Proposition \ref{fix-points} and \ref{prop:fix-points} already stated in the literature}\label{proof-old}

\begin{lem}[Particular case of Theorem 2.1 \cite{sergeantperthuis2022regionalized}]\label{caracterization-critical-points}

Let $F$ be a presheaf from a finite poset $\A$ to the category of finite sets $\textbf{FinSet}$. A point $v\in \lim \widetilde{F}$ is a critical point of $F_{\text{Bethe}}$ if and only if

\begin{equation}\label{thm-local-optimization-critical-points}
\big[\mu_{F^\ast} d_v \text{FE}\big] \big|_{\lim \widetilde{F}} = 0,
\end{equation}

where $|_{\lim \widetilde{F}}$ means restricted to the vector space $\lim \widetilde{F}$.

The previous Equation \eqref{thm-local-optimization-critical-points} is equivalent to the existence of $l\in F_\rightarrow$ such that

\begin{equation}\label{critical-points}
\mu_{F^\ast} d_v \text{FE} = d_{F^\ast}(l).
\end{equation}

\end{lem}
\begin{proof}
One shows that for any $v\in \bigoplus_{a\in \A}\widetilde{F}_a$,  

\[
\text{d}_v F_{\text{Bethe}}\big|_{\lim \widetilde{F}} = \big[\mu_{F^\ast} d_v \text{FE}\big] \big|_{\lim \widetilde{F}},
\]

where  

\[
\text{FE}_a(h_a) = \sum_{x_a} h_a(x_a) H_a(x_a) + \sum_{x_a} h_a(x_a) \ln h_a(x_a),
\]

for $h \in \R^{E_a}$. Its differential is given by Equation \eqref{proof:differential}. Furthermore, a parametrization of $\lim \widetilde{F}$ is given by $\im \text{d}_{F^\ast}$.
\end{proof}

As $\zeta_{F^\ast}$ is the inverse of $\mu_{F^\ast}$, Equation \eqref{critical-points} is equivalent to the existence of $l\in F_\rightarrow$ and $v\in \lim \widetilde{F}$ such that  

\[
d_v \text{FE} = \zeta_{F^\ast} d_{F^\ast}(l).
\]

\subsection{Proof of Proposition \ref{fix-points}}

The belief propagation algorithm Equation \ref{BP-update-rule} can be rewritten as, 

\[
\ln m^{t+1}_{a\to b} = \ln m^{t}_{a\to b} + \ln \widetilde{F}^a_b 
\left( g_H \circ \zeta_{F^{\dagger}} \circ d_{F^\dagger} (\ln m^{t})_a \right)  
- \ln \left( g_H \circ \zeta_{F^{\dagger}} \circ \text{d}_{F^\dagger} (\ln m^{t})_b \right)
\]

Let us denote $\Delta \text{BP}(\ln m)_{a\to b}= \ln \widetilde{F}^a_b 
\left( g_H \circ \zeta_{F^{\dagger}} \circ d_{F^\dagger} (\ln m)_a \right)  
- \ln \left( g_H \circ \zeta_{F^{\dagger}} \circ \text{d}_{F^\dagger} (\ln m)_b\right)$. 

Consider the collection $\left(\ln m_{a\to b}+C_{a\to b}; a,b\in \A: b\leq a\right)$, where $C_{a\to b}$ is a constant, i.e. it does not depend on $x_b\in F_b$. Then there is a collection of constants $(C^{'}_{a\to b}; a,b\in \A: b\leq a)$ such that $\Delta \text{BP}(\ln m_{a\to b}+C_{a\to b})= \Delta \text{BP}( \ln m_{a\to b})+ C^{'}_{a\to b}$. Therefore $\text{BP}$ is a algorithm that preserves the equivalence classes $\{m+C\}$, we will denote that class as $[m]$ and $\overline{\Delta \text{BP}} ([m])= [\Delta \text{BP}(m)]$.

Let $[\ln m]^{\ast}$ be a fix point of the Belief Propagation algorithm, and let $\ln m^{r}$ be a representant of $[\ln m]^{\ast}$ i.e. $[\ln m]^{\ast}= [\ln m^{r}]$. There is a a collection of constants $(C_{a\to b}; a,b\in \A: b\leq a)$, 

$$C_{a\to b}= \ln \widetilde{F}^a_b 
\left( g_H \circ \zeta_{F^{\dagger}} \circ \text{d}_{F^\dagger} (\ln m^{r})_a\right)-\ln g_H \circ \zeta_{F^{\dagger}} \circ \text{d}_{F^\dagger} (\ln m^{r})_b $$

Therefore, 

$$\widetilde{F}^a_b 
\left( g_H \circ \zeta_{F^{\dagger}}\circ \text{d}_{F^\dagger} (\ln m^{r})_a\right)= e^{C_{a\to b}}g_H \circ \zeta_{F^{\dagger}} \circ \text{d}_{F^\dagger} (\ln m^{r})_b$$

The collection of constant function $\R\hookrightarrow \widetilde{F}_a$ in each space $F_a$ with $a\in \A$ forms a subobject of $F$, therefore $[\widetilde{F}](a)= \widetilde{F}_a/ \mathbb{R}$, with associated quotient maps $[\widetilde{F}]^a_b$ is a presheave from $\A$ to finite vector spaces. A probability distribution $Q\in \widetilde{F}_a$ is in correspondance with $[Q]\in  [\widetilde{F}](a)$, the section $b\in \widetilde{F}_a\to b/\sum b$ is a section which in turn provides the inverse map when passing to the quotient. The previous equation implies that, 

$$[\widetilde{F}^a_b 
\left( g_H \circ \zeta_{F^{\dagger}}\circ \text{d}_{F^\dagger} (\ln m^{r})_a\right)]=[ g_H \circ \zeta_{F^{\dagger}} \circ \text{d}_{F^\dagger} (\ln m^{r})_b]$$

Therefore, we have just shown that \( \overline{\Delta \text{BP}}([m^r]) = 0 \) implies that  
\[
[g_H \circ \zeta_{F^{\dagger}} \circ \text{d}_{F^\dagger} ]([\ln m]^*) \in \lim [\widetilde{F}].
\]  
Define \( b_a = (g_H \circ \zeta_{F^{\dagger}} \circ \text{d}_{F^\dagger}([\ln m]^*))_a\).  
This statement can be reformulated as saying that  
\[
\left( \frac{b_a}{\sum_{x_a \in F_a} b_a}, a \in \mathcal{A} \right) \in \lim (\mathbb{P} F).
\]

The previous equation implies that $\left( \frac{b_a}{\sum_{x_a \in F_a} b_a}, a \in \mathcal{A} \right) $ is a critical point of the Bethe Free Energy.

Conversely, if one starts with a critical point of the Bethe Free Energy, \( p_a \in \mathbb{P}(F_a), \, a \in \mathcal{A} \), then there exists a collection \( (m^{\ast}_{a \to b}; \, a, b \in \mathcal{A}, \, b \leq a) \) such that  

\[
p = g_H \circ \zeta_{F^\dagger} \circ d_{F^{\dagger}} (m^{\ast}).
\]

Since \( p \in \lim \mathbb{P}F \), we have  

\[
\ln \widetilde{F}^a_b(p_a) - \ln p_b = 0
\]

and  

\[
\text{BP}(m^{\ast}) = m^{\ast}.
\]

\qed

\subsection{Proof of Proposition \ref{prop:fix-points}}
We will follow the proof of Theorem 2.1 \cite{sergeantperthuis2022regionalized}. Let $m^{\ast}$ be a fixed point of $\text{MP}$, i.e., $\Delta \text{MP}(m^\ast)= 0$. This implies that  
$$g_H\circ \zeta_{F^{\dagger}} \circ \text{d}{F^\dagger} (\ln m^{\ast})\in \lim \widetilde{F},$$  
and that the associated normalized belief is in $\lim\mathbb{P}F$. By Lemma~\ref{caracterization-critical-points}, therefore, $m^{\ast}$ is a critical point of the Bethe Free Energy.
\qed

\section{Link between belief propagation and message passage algorithms}\label{comparison-algorithm}

When we study the functoriality of message passing algorithms, we consider the following operator:  

$$\Delta \text{MP}(l)= \delta_F \circ g_H\circ \zeta_{F^\dagger}\circ d_{F^\dagger}(l)$$

The messages $l$ of $\text{MP}$ are the $\ln m$ of $\text{BP}$, with for any $a,b\in \A$ such that $b\leq a$ and any $x_b\in F_b$, $m_{a\to b}(x_b)=e^{l_{a\to b}(x_b)}$. Let us now show that for any fixed point $l^{\ast}$ of $\text{MP}$, $[m^{\ast}]$ is a fixed point of $\text{BP}$. We will then show that for any fixed point $[m]^{\ast}$ of $\text{BP}$, there exists $m^{\ast} \in [m]^{\ast}$ such that $\ln m^\ast$ is a fixed point of $\text{MP}$. Showing this justifies why the study of the critical points of the Belief Propagation algorithm is equivalent to the study of the $\text{MP}$ algorithm considered in \cite{sergeantperthuis2022regionalized} and in this paper.

Let $\text{MP}(l^\ast)= l^\ast$ then $\Delta \text{MP} (l^\ast)=0$ and so for any $b\leq a$, $F^a_b\left(g\circ \zeta_{F^\dagger} \circ d_{F^{\dagger}}(l)_a\right)=g\circ \zeta_{F^\dagger} \circ d_{F^{\dagger}}(l)_b$. Therefore, in particular $\Delta \text{BP}(\ln m^{\ast})=0$ which implies that $m^\ast$ is a fix point of $ \text{BP}$.\\

For any message \( \ln m_{a\to b} \in \widetilde{F}_b \), define \( N_{a\to b}(\ln m_{a\to b}) = \ln m_{a\to b} - \ln \sum_{x_b} e^{\ln m_{a\to b}} \). Then, \( N_{a\to b}(\ln m_{a\to b} + C) = N_{a\to b}(\ln m_{a\to b}) \). Let us denote \( N(\ln m) = (N_{a\to b}; a \to b) \).  

Consider \( [m]^{\ast} \) a fixed point of \( \text{BP} \). Then, for any representative \( m^{r} \) such that \( [m^{r}] = [m]^{\ast} \), we have \( \Delta \text{MP}(N(\ln m^{r})) = 0 \); therefore, \( N(\ln m^{r}) \) is a fixed point of \( \text{MP} \).

\section{Proof of main theorems}\label{main-proofs}

\subsection{Proof of Theorem \ref{main:thm1}}

Now it is easy to show that $\phi_1\circ\delta_F= \delta_G\circ \phi$ therefore by dualizing (as $d_F:=(\delta_F)^*= d_{F^*}$) on ehas that $d_F\circ \phi^{*}_1= \phi^{*}  \circ d_G $. \\

Let us now show that $\zeta$ also commutes to natural transformations for functors that have the same domain: 

$$\zeta_{F^*}\circ \phi^*(v)(a)=\sum_{b\leq a} {F^a_b}^{*}\phi^*(v)(b)$$

$$\sum_{b\leq a} {F^a_b}^{*}\phi^*_b v_b=\sum_{b\leq a} {\phi_b \circ F^a_b}^{*}v_b$$ 
and $\sum_{b\leq a} {\phi_b \circ F^a_b}^{*}v_b= \sum_{b\leq a} {G^a_b \circ \phi_a}^{*}v_b=\sum_{b\leq a} \phi_a^{*}\circ {G^a_b}^{*} v_b $ and $\phi_a$ is linear so, $\sum_{b\leq a} \phi_a^{*}\circ {G^a_b}^{*} v_b = \phi_a^{*}\sum_{b\leq a}  {G^a_b}^{*} v_b$. So $\zeta_{F^*}\circ \phi^*= \phi^*\circ \zeta_{G^*}$. \\

$g_H: \bigoplus_{a \in \mathcal{A}} \widetilde{F}_a \to \bigoplus_{a \in \mathcal{A}} \widetilde{F}_a
$ is defined by the following relationship: 

\begin{equation}
d_xf_a=y \iff x=g_a(y)
\end{equation}

where $f_a(q_a)=\sum q_a H_a - S(q_a)= \sum_{\omega} q_a(\omega) H_a(\omega) +\sum_{\omega} q_a(\omega) \ln q_a(\omega)$. The differential gives for the choice of scalar product 

$$y_a= H_a +\ln q_a +1$$

What we want to show is that there is $\widetilde{H}:=(\widetilde{H}_a: G_a\to \mathbb{R},a\in \mathcal{A})$ such that $\phi\circ g_{H}\circ \phi^\dagger =  g_{\widetilde{H}}$. We must show it ponctually as we do not requiere any compatibility condition on the $H_a$: i.e. we need to show that, $\phi\circ g_{H_a}\circ \phi^\dagger= g_{\widetilde{H_a}}$.
Let us first remark that for any $p_a\in \widetilde{F}_a$, $\phi^{\dagger}_a(p_a)=p_a\circ \phi_a$. Then,

\[
\phi_a\circ g_{H_a}\circ \phi^{\dagger}_a(p)(y_a) = \sum_{x_a:\phi_a(x_a)=y_a} e^{-H_a(x_a) + p\circ \phi(x_a) -1} 
= e^{p(y_a)-1} \sum_{x_a:\phi_a(x_a)=y_a} e^{-H_a(x_a)}.
\]
\[
\phi_a\circ g_{H_a}\circ \phi^{\dagger}_a(p)(y_a) = \mathbf{1}[y_a\in \operatorname{im} \phi_a] e^{p(y_a)-1-\widetilde{H}_a}.
\]
By adopting the convention that \( e^{-\infty}=0 \), we write:
\[
\phi_a\circ g_{H_a}\circ \phi^{\dagger}_a(p)(y_a) = e^{p(y_a)-1-\widetilde{H}_a}= g_{\widetilde{H}_a}.
\]

W proved that if there is a natural transformation of measurable map between $\phi: F\to G$ then, 

$$\Delta MP_{G}= \phi\circ \Delta MP_{F}\circ \phi^{\dagger}$$

and so one deduces from the fact that $g_{\widetilde{H}_a}^{\langle\cdot ,\cdot\rangle}= \phi_a\circ g_{\widetilde{H}_a}^{\langle\cdot ,\cdot\rangle}\circ \phi^{\dagger}_{a,\langle\cdot,\cdot\rangle}$ that $$\Delta MP_{G}^{\langle \cdot,\cdot\rangle}= \phi\circ \Delta MP_{F}^{\langle \cdot,\cdot\rangle}\circ \phi^{\dagger}$$.

\subsection{Proof of Theorem \ref{fonctoriality-increment}}

Consider $\phi:F\to G$ a natural transformation, then $\Delta \text{MP}(\phi)=(\phi,\psi)$ and by Theorem\ref{main:thm1}, 

$$\Delta \text{MP}(G)(l,\psi(H))= \phi\circ \Delta \text{MP}_H(\phi^\dagger (l))$$

Furthermore, for any \( a \in \mathcal{A} \) and \( y_a \in G_a \), the mapping \( H \mapsto \psi(H)_a(y_a) \) is convex and therefore continuous. Since \( \mathcal{A} \) is finite and each \( G_a \) is finite, \( \psi \) is a continuous map. We just showed that $\Delta \text{MP}(\phi): \Delta \text{MP}(F)\to \Delta \text{MP} (G)$ is a morphism of $\textbf{MPA}$. 

Let us now show that \( \Delta \text{MP} \) is functorial. Let \( F, G, L \) be three presheaves of \( \hat{\mathcal{A}}_f \), and let \( \phi: F \to G \), \( \phi_1: G \to L \) be two natural transformations. Then,

\[
\forall a \in \mathcal{A}, \, z_a \in L_a, \quad \Delta \text{MP}(\phi_1) \circ \Delta \text{MP}(\phi)(z_a) = -\ln \sum_{y_a \in G_a \,:\, \phi_{1,a}(y_a) = z_a} e^{-\psi_a(y_a)}
\]

\[
\begin{aligned}
-\ln \sum_{\substack{y_a \in G_a \\ \phi_{1,a}(y_a) = z_a}} e^{-\psi_a(y_a)}  
&= -\ln \sum_{\substack{y_a \in G_a \\ \phi_{1,a}(y_a) = z_a}} 
\sum_{\substack{x_a \in F_a \\ \phi_{a}(x_a) = y_a}} e^{-H_a(x_a)} \\
&= -\ln \sum_{\substack{x_a \in F_a \\ \phi_{1} \circ \phi_a(x_a) = z_a}} e^{-H_a(x_a)}
\end{aligned}
\]

 Therefore $\Delta \text{MP}$ is a functor. 

 \qed

\subsection{Proof of Theorem \ref{label:thm-3}}
Let $F,G$ be two presheaves of $\hat{A}_f$, let $\phi: F\to G$ be a monomorphism, i.e. a natural transformation where each map is surjective. Let $\text{MP}(\phi)=(\phi,\psi)$ with, 
\[
\forall a \in \mathcal{A}, \, y_a \in G_a, \quad \psi(H)_a(y_a) = -\ln \sum_{x_a \in F_a : \phi_a(x_a) = y_a} e^{-H_a(x_a)}.
\]

By Theorem \ref{fonctoriality-increment}, $\phi\circ \Delta \text{MP}_{F,H}\circ \phi^\dagger= \Delta \text{MP}_{G,\psi(H)}$.

$$\phi \circ \text{MP}_{F,H} \circ \phi^{\dagger}(l) = \phi \circ \phi^\dagger(l) + \phi \circ \text{MP}_{F,H} \circ \phi^{\dagger}(l)
$$

For the choice of scalar product \( \langle \cdot, \cdot \rangle^\phi \), by Lemma \ref{lem:isometry}, \( \phi \circ \phi^\dagger = \text{id} \), therefore, 
$$
\phi \circ \text{MP}_{F,H} \circ \phi^{\dagger}(l) = l + \Delta \text{MP}_{F,\psi(H)}(l) = \text{MP}_{G,\psi(H)}(l),
$$

and therefore \( \text{MP}(\phi): \text{MP}(F) \to \text{MP}(G) \) is a morphism of \( \text{MPA} \); functoriality is a consequence of Theorem \ref{fonctoriality-increment}, as \( \Delta \text{MP}(\phi) = \text{MP}(\phi) \).

\end{document}